\documentclass[11pt]{article}

\usepackage{vien}

\usepackage[utf8]{inputenc} % allow utf-8 input
\usepackage[T1]{fontenc}    % use 8-bit T1 fonts
\usepackage{mathtools,subfigure,bm,epstopdf,url}
\usepackage{algorithm}    
\usepackage{enumitem}
\usepackage[final]{changes}

\usepackage[top=2.54cm,right=3cm,left=3cm,bottom=2.54cm]{geometry}

\newcommand{\intr}{\mathop{\mathrm{int}}}

\newcommand{\genf}{\varphi}
\newcommand{\breg}[1]{D_{\genf}\left(#1\right)}

\newcommand{\ones}{\mathbf{1}}

\newcommand{\ya}{y^{\alpha}}
\newcommand{\gam}[2]{\gamma_{#1}^{{#2}}}
\newcommand{\opt}{^\star}

\title{Anderson Acceleration of Proximal~Gradient~Methods}

\author{	
	Vien V.~Mai\footnotemark[1]
	\and  
	Mikael Johansson{\thanks{Division of Decision and Control Systems, School of Electrical Engineering and Computer Science, KTH~Royal Institute of Technology, SE-100 44  Stockholm, Sweden. Emails: \tt\small\{maivv, mikaelj\}@kth.se.}}       
}

\begin{document}
\maketitle
\begin{abstract}
Anderson acceleration is a well-established and simple technique for speeding up fixed-point computations with countless applications. This work introduces novel methods for adapting Anderson acceleration to proximal gradient algorithms. Under some technical conditions, we extend existing local convergence results of Anderson acceleration for smooth fixed-point mappings to the proposed non-smooth setting. We also prove analytically that it is in general, impossible to guarantee global convergence of native Anderson acceleration. We therefore propose a simple scheme for stabilization that combines the global worst-case guarantees of proximal gradient methods with the local adaptation and practical speed-up of Anderson acceleration. Finally, we  provide the first applications of Anderson acceleration to non-Euclidean geometry.
\end{abstract}

\section{Introduction}
The last few decades have witnessed significant advances in the
theory and practice of convex optimization based on first-order  information \cite{Nes04, Bec17}. The
worst-case oracle complexity has been established for many function classes \cite{NY83} and  algorithms with matching worst-case performance have been developed. 
However, these methods are only optimal in a worst-case (resisting oracle) sense, and are developed under the assumption that global function properties are known and constant. In practice, however, such constants are almost never known a priori. Moreover, their local values, which determine the actual practical performance, may be very different from their conservative global bounds and often change as the iterates approach optimum. It is also observed that acceleration methods such as Nesterov's accelerated gradient are very sensitive to misspecified parameters; slightly over- or under-estimating the strong convexity constant can have a severe effect on the overall performance of the algorithm \cite{OC15}.
Thus, strong practical performance of optimization algorithms requires \emph{local adaption and acceleration}. Efficient line-search procedures \cite{Nes13}, adaptive restart techniques \cite{OC15} and nonlinear acceleration schemes \cite{SdB16} are therefore now receiving an increasing attention.

Extrapolation techniques have a long history in numerical analysis (see, e.g.,~\cite{Sid17, BRS18}). Recently, its idea has resurfaced in the first-order optimization literature \cite{SdB16, ZOB18, MJ19, FZB19, PL19}. Unlike momentum acceleration methods such as Polyak's heavy ball~\cite{Pol64} and Nesterov's fast gradient~\cite{Nes04}, which require knowledge of problem parameters, classical extrapolation techniques for vector sequences such as minimal polynomial extrapolation~\cite{SFS87}, reduced rank extrapolation~\cite{Edd79}, vector epsilon algorithm~\cite{Wyn62}, and Anderson acceleration~\cite{And65} estimate the solution directly from the available \added{iterate} sequence. 
These methods enjoy favorable theoretical properties of Krylov subspace methods on quadratic problems and often perform equally well in practice on non-quadratic problems.

\subsection{Related Work}
Anderson acceleration \added{(AA)} was proposed in the 1960's to expedite solution times for nonlinear integral equations \cite{And65}. The technique has then been generalized to general fixed-point equations and found  countless applications in diverse fields such as computational chemistry, physics, material science, etc.~\cite{Pul80, Eye96, WN11}.
However, AA and optimization algorithms have \replaced{been}{so far} developed quite independently and only limited connections were discovered and studied~\cite{Eye96, FS07}. Very recently, the technique has started to gain a significant interest in the optimization community (see, e.g.,~\cite{SdB16, SdB17, BSd18, ZOB18, FZB19, PL19}). Specifically, a series of papers~\cite{SdB16, SdB17, BSd18}  adapt AA to accelerate several classical algorithms for unconstrained optimization; \cite{ZOB18} studies a variant of AA for non-expansive operators; \cite{FZB19} proposes an application of AA to Douglas-Rachford splitting; and \cite{PL19} uses AA to improve the performance of the ADMM method. There is also an emerging literature on applications of AA in machine learning~\cite{HS16, MSG18, GS18, PDZ18}. 

Although some initial success has been obtained for adapting AA to optimization algorithms, current research has mainly focused on unconstrained or linearly constrained minimization (e.g., \cite{SdB16, FZB19}). For non-smooth composite problems, \emph{asymptotic} convergence results of AA are often achieved by additional safeguarding strategies \cite{ZOB18}, without which even local convergence guarantees have not been available. This is because AA relies on linearization (and hence often requires differentiability) of the associated mapping around its fixed-point, which is hard to adapt to non-smooth  optimization. Our aim with this paper is to address these limitations.  To this end, we make the following contributions:
\begin{description}[leftmargin=0.45cm]
\item [1.] We propose a simple and efficient AA scheme for the classical proximal gradient algorithm (PGA) and, under mild technical conditions, establish local convergence.

\item [2.] Local convergence properties of \emph{native} AA have been studied in various settings~\cite{TK15,SdB16, PR19, LL18, MJ19}. However, whether native AA converges globally still remains largely unknown (cf.~\cite{FZB19}). Here, we show a negative answer to this question. More specifically, we construct an unconstrained strongly convex problem for which we can prove analytically that AA fails to converge.  
We therefore stabilize the proposed method by a simple guard step that preserves the global worst-case convergence guarantees of PGA without sacrificing the local adaption and acceleration abilities of AA.

\item [3.] We adapt AA to the Bregman proximal gradient (BPG) family, where the mirror descent \cite{NY83} and NoLips \cite{BBT16} methods are special instances. The method respects the structure of the BPG family and admits a simple and elegant interpretation. To the best of our knowledge, these are the first applications of AA to non-Euclidean geometry. 

\item[4.] We perform substantial experiments on several important classes of constrained optimization problems and demonstrate consistent and dramatic speedups on real-world data-sets. 

\end{description}

\subsection{Notation}
We denote by $\R_+$ the set of nonnegative real numbers. For a set $\mathcal{X}$, $\overline{\mathcal{X}}$  and $\intr \mathcal{X}$ denote  its closure and  interior, respectively. The notation $\norm{\cdot}$ refers to a general norm, and $\norm{\cdot}_2$ is the Euclidean norm. The all-ones vector is denoted by $\ones$. Finally, the vector quantity $x = o(t)$ with $t>0$ means that $\ltwo{x}/t\to 0$ as $t\to 0$.

\section{Anderson acceleration}
Let $g: \R^n \to \R^n$ be a mapping and consider the problem of finding a fixed-point of $g$: 
\begin{align*}
	\mbox{Find} \,\, x \in \R^n \,\, \mbox{such that}\,\, x = g(x). 
\end{align*} 
%Let $\{x_k\}$ be a  sequence of iterates generated by the fixed-point iteration: 
%\begin{align*}
%	x_{k+1}=g(x_k), \quad k=0,1,2,\ldots.
%\end{align*}
In contrast to the fixed-point iteration $y_{k+1}=g(y_k)$,  which only uses the last iterate to generate a new estimate, AA tries to make better use of past information. Concretely, let $\{x_i\}_{i=0}^{k}$ be the sequence of iterates generated by AA up to iteration $k$. Here, we refer the term $r_k := g(x_k)-x_k$ as the residual in the $k$th iteration. Then, to form $x_{k+1}$, it searches for a point that has smallest residual within the subspace spanned by the $m+1$ most recent iterates. In other words, if we let $\bar{x}_k=\sum_{i=0}^m \alpha_i^k x_{k-i}$, AA seeks to find a vector of coefficients $\alpha^k=[\alpha_0^k,\ldots,\alpha_{m}^k]^\top$  such that
\begin{align}\label{eq:ideal:coeff}
	\alpha^k
	=
		\argmin_{\alpha: \alpha^\top\ones=1} \big\| g(\bar{x}_k) -\bar{x}_k\big\|.
\end{align}
However, since~\eqref{eq:ideal:coeff} can be hard to solve for a general nonlinear mapping $g$, AA uses
%solves instead 
\begin{align}\label{eq:coeff}
	\alpha^k 
	=
		\argmin_{\alpha: \alpha^\top\ones=1}
		\Big\|{\sum_{i=0}^m \alpha_i g(x_{k-i})-\sum_{i=0}^m \alpha_i x_{k-i}}\Big\|.
%	=
%		\argmin_{\alpha: \alpha^\top\ones=1}
%		\Big\|{\sum_{i=0}^k \alpha_i r_i}\Big\|.
\end{align}
It is clear that Problems~\eqref{eq:ideal:coeff} and~\eqref{eq:coeff} are  equivalent if $g$ is an affine mapping. Let $R_k=[r_k,\ldots,r_{k-m}]$ be the residual matrix at the $k$th iteration, Problem~\eqref{eq:coeff} can then be written as
\begin{align}\label{eq:ls}
  \alpha^k= \argmin_{\alpha^\top\ones=1}\norm{R_k\alpha}.
\end{align}
With $\alpha^k$ computed, the next iterate of AA is then generated by 
\begin{align}
	x_{k+1} = \sum_{i=0}^{m} \alpha_i^k g\left(x_{k-i}\right),
\end{align}
which in the affine case, is equivalent to applying the operator $g$ to $\bar{x}_k$. When $m=0$, AA reduces to the fixed-point iteration.

\begin{algorithm}[!h]
	\caption{Anderson Acceleration}
	\begin{algorithmic}[1]\label{alg:aa}
		\REQUIRE $x_0$, $m\geq 0$, $g(\cdot)$			
		\STATE $x_1\gets g(x_0)$
		\FOR{$k=1,\ldots, K-1$}
		\STATE 	$m_k \gets \min(m,k)$
		\STATE  $R_k \gets [r_k, \ldots, r_{k-m_k}]$, where $r_i = g(x_i)-x_i$
		\STATE $ \alpha^k	\gets \argmin_{\alpha^\top\ones=1}\norm{R_k \alpha}$
		\STATE $x_{k+1} \gets  \sum_{i=0}^{m_k} \alpha_i^k g(x_{k-i})$
		\ENDFOR
		\ENSURE $x_K$
	\end{algorithmic}
\end{algorithm}

One of the reason that AA is so popular in engineering and scientific applications is that it can speed-up convergence with almost no additional tuning parameters and the  extrapolation coefficients can be computed very efficiently. When the Euclidean norm is considered, Problem \eqref{eq:ls} is a simple least-squares, which admits a closed-form solution given by
\begin{align}
	\alpha^k = \frac{(R_k^\top R_k )^{-1}\ones}{\ones^\top (R_k^\top R_k )^{-1} \ones}.
\end{align}
This can be solved by first solving the $m \times m$ normal equations $R_k^\top R_k x=\ones$ and then normalizing the result to obtain $\alpha^k=x/(\ones^\top x)$ \cite{SdB16}. Indeed, the computations can be done even more efficiently using QR decomposition. 
When passing from $R_{k-1}$ to $R_{k}$, only the last column of $R_{k-1}$ is removed and a new column is added. Thus, the corresponding $Q$ and $R$ matrices can be easily updated and the total cost is at most $O\left(m^2 + mn\right)$ \cite{HS16}. Since $m$ is typically between $1$ and $10$ in practice, this additional cost is negligible compared to the cost of evaluating $g$. 
%\begin{figure*}[h!]
%	\centering 
%	{\includegraphics[width=0.7\textwidth]{Figures/quadractic_cond_1e3.png}}
%	{\includegraphics[width=0.7\textwidth]{Figures/quadractic_cond_1e4.png}}
%	~\caption{Performance of difference first-order optimization algorithms on quadratic convex problems. Here, $A\in\R^{100\times 100}$ is a symmetric positive semidefinite matrix with 25 nonzero eigenvalues. On the top $\lambda_1(A)/\lambda_{25}(A)=10^3$ and on the bottom $\lambda_1(A)/\lambda_{25}(A)=10^4$.}~\label{fig:quadratic} 
%\end{figure*}

\begin{figure}[h!]
	\centering 
	\begin{minipage}{0.5\textwidth}
		\centering 
		{\includegraphics[width=\textwidth]{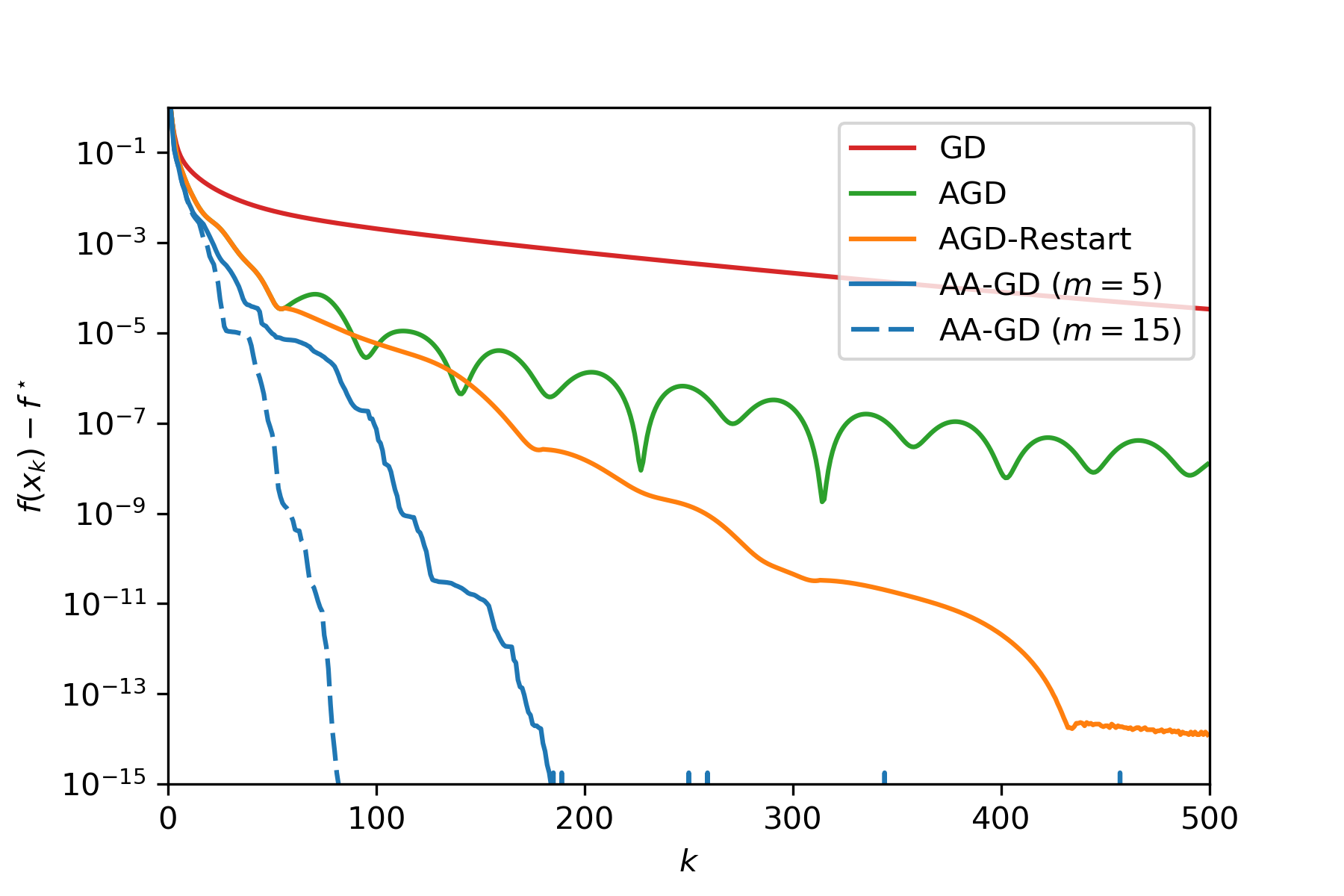}}	
	\end{minipage}
	\hskip -0.1in
	\begin{minipage}{0.5\textwidth}
		\centering 
		{\includegraphics[width=\textwidth]{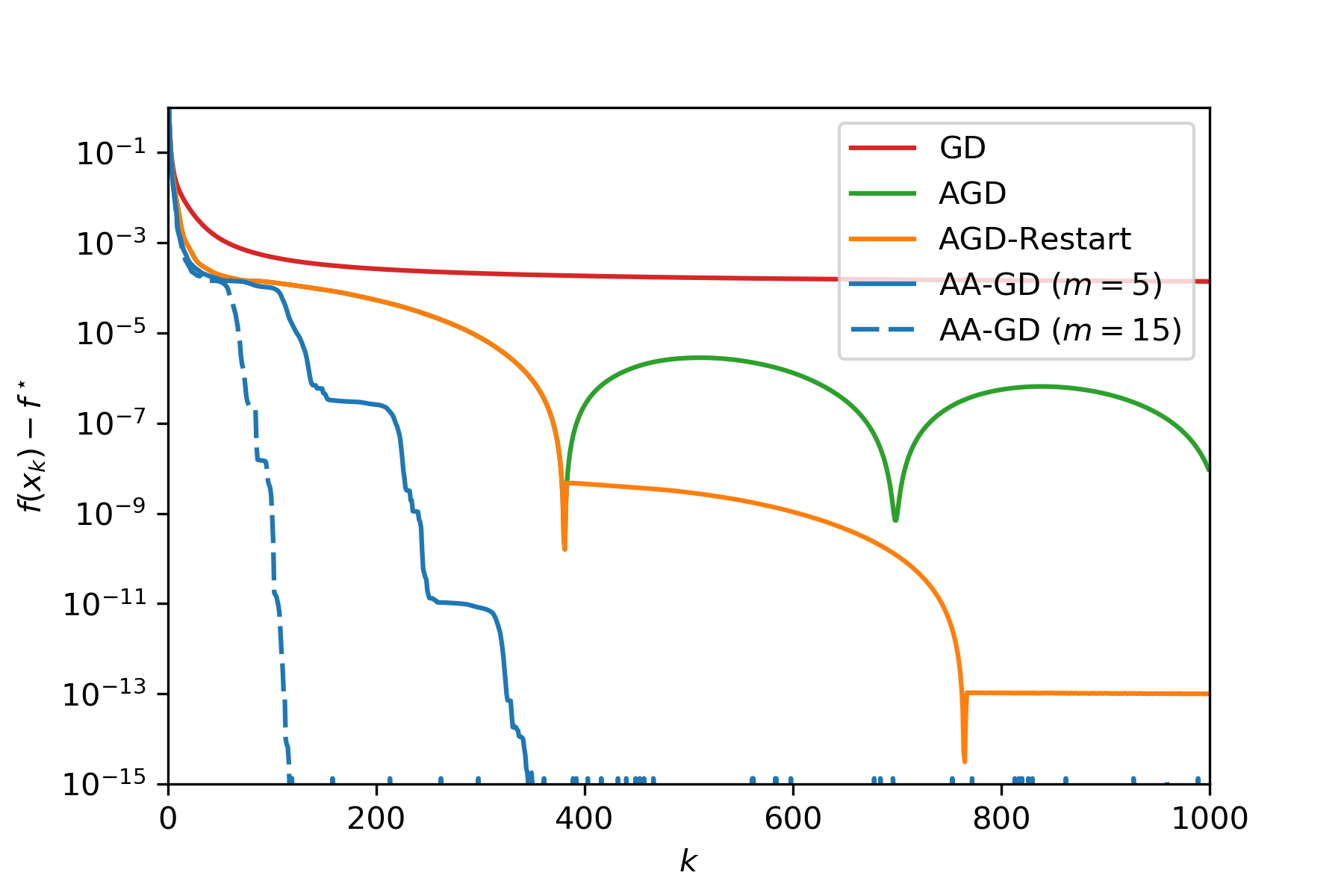}}
	\end{minipage}
	\caption{Quadratic convex problems: Left: $\lambda_1(A)/\lambda_{25}(A)=10^3$. Right: $\lambda_1(A)/\lambda_{25}(A)=10^4$.}~\label{fig:quadratic} 
\end{figure}

\subsection{Anderson acceleration for optimization algorithms}
Since many optimization algorithms can be written as fixed-point iterations, they can be accelerated by the memory-efficient, line search-free AA method with almost no extra cost. 
For example, the classical gradient descent (GD) method for minimizing a smooth convex function $f$ defined by
\begin{align*}
x_{k+1} = x_k - \gamma\grad{f(x_k)},
\end{align*}
is equivalent to the fixed-point iteration applied to $g(x)=x - \gamma\grad{f(x)}$. Clearly, a fixed-point of $g$ corresponds to an optimum of $f$. 
The intuition behind AA for GD is that smooth functions are well approximated by quadratic ones around their (unconstrained) optimum, so their gradients and hence $g$ are linear. In such regimes, AA enjoys several nice properties of Krylov subspace methods.
Specifically, consider a convex quadratic minimization problem
\begin{align}
	\underset{x\in \R^n}{\mbox{minimize}} \,\, \frac{1}{2} x^\top A x - b^\top x,
\end{align}
where $A \in \R^{n \times n}$ is a symmetric positive semidefinite matrix and $b \in \R^n$. It has been shown in \cite{WN11, PE13} that AA with full information (i.e. setting $m=\infty$ in Step~3 of Algorithm~\ref{alg:aa}) is essentially equivalent to GMRES \cite{SS86}.  Therefore, AA admits the convergence rate~\cite{Nem95, Gre97} 
\begin{align}\label{eq:opt:rate}
	\norm{x_k - x^\star}_2^2 
	\leq 
		O\left(\min\left\{1/k^2,e^{-k/\sqrt{\kappa}} \right\}\right)\norm{x_0-x^\star}_2^2,
\end{align} 
where $\kappa=\lambda_1(A)/\lambda_n(A)$ is the condition number.
This rate shows a very strong adaptation ability and is attained without any knowledge of the problem at hand, a remarkable property of Krylov subspace methods. 
In contrast, Nesterov’s accelerated gradient method (AGD) \cite{Nes04} can only achieve this rate if $\lambda_1(A)$ and $\lambda_n(A)$ are known. 
%When $m$ is finite, it can be shown that AA will not do worse than GD for any $m$. More concretely, for any $m \in \{0,1,\ldots\}$, AA achieves the worst-case guarantee
%\begin{align*}
%	\norm{x_k - x^\star}_2^2 
%	\leq 
%		O\left(e^{-k/{\kappa}}\right)\norm{x_0-x^\star}_2^2,
%\end{align*}
%which coincides with that of GD. The above observations imply that 
%AA interpolates the two extreme regimes where $m=0$ and $m=\infty$. 

%In practice, significant speed-ups and strong adaptation are often observed even with very small $m$.
%As an example, Figure~\ref{fig:quadratic} shows the performance of different algorithms applied to minimize a quadratic convex function in $n=100$ dimensions with $25$ nonzero eigenvalues such that $\lambda_1(A)/\lambda_{25}(A)$ equals $10^3$  and $10^4$ for the top and the bottom plots, respectively. We compared AA-GD with GD, AGD, and the adaptive restart scheme (AGD-Restart) in \cite{OC15}. It should be noted that just like AA, the main objective of the AGD-Restart scheme is to achieve local adaptation. We can see that local adaptation and acceleration can dramatically improve the performance of an optimization algorithm. It is evident that AA initially converges at the same rate as AGD ($1/k^2$) and eventually switches to linear convergence, as suggested in \eqref{eq:opt:rate}, even with a very small value of $m$ and a non-strongly convex objective. 

In practice, significant speed-ups and strong adaptation are often observed even with very small $m$.
As an example, Figure~\ref{fig:quadratic} shows the performance of different algorithms applied to minimize a quadratic convex function in $n=100$ dimensions with $25$ nonzero eigenvalues.
We compared AA-GD with GD, AGD, and the adaptive restart scheme (AGD-Restart) in \cite{OC15}. It should be noted that just like AA, the main objective of the AGD-Restart scheme is to achieve local adaptation. We can see that local adaptation and acceleration can dramatically improve the performance of an optimization algorithm. It is evident that AA initially converges at the same rate as AGD ($1/k^2$) and eventually switches to linear convergence, as suggested in \eqref{eq:opt:rate}, even with a very small value of $m$ and on an objective function which is not strongly convex.

%\noindent\textbf{Extension to smooth convex minimization:}
If the function being minimized has a positive definite Hessian at the optimum, then near the solution it can be well approximated by a quadratic model
\begin{align*}
	f(x) \approx f(x^\star) + (x-x^\star)^\top \nabla^2f(x^\star)(x-x^\star).
\end{align*}
Note that the matrix $\nabla^2f(x^\star)$ may have smallest eigenvalue $\lambda_{\mathrm{min}}$ strictly greater than the \emph{global} strong convexity constant $\mu$. Thus, once we enter this regime, we may be able to achieve all the nice features of AA on quadratic problems discussed in the previous paragraphs. 

\subsection{Anderson acceleration as a multi-step method}
It is known that AA is related to several iterative schemes such as multisecant quasi-Newton methods~\cite{Eye96,FS07,WN11}. Here, we point out some connections between AA-GD and multi-step methods in optimization. To do so, let $\gam{i}{k}:=\sum_{j=i}^{m_k} \alpha_j^k$, $i \in \{1,\ldots,m_k\}$ and define 
$\ya_k :=  \sum_{i=1}^{m_k} \alpha_i^k x_{k-i}$. AA-GD can \added{then} be written as
\begin{align}\label{eq:ya}
	\ya_k = x_k - \sum_{i=1}^{m_k}\gam{i}{k}\left(x_{k-i+1}-x_{k-i}\right)
	\quad \mbox{and} \quad
	x_{k+1} = \ya_k - \gamma \sum_{i=0}^{m_k} \alpha_i^k \grad f(x_{k-i}).
\end{align}
Recall that Nesterov's accelerated gradient method (AGD)~\cite{Nes04} can be written as
\begin{align*}
	y_k = x_k + \beta_k(x_k-x_{k-1})
	\quad \mbox{and} \quad
	x_{k+1} = y_k-\gamma\grad{f}(y_k),
\end{align*} 
while Polyak's  Heavy ball (HB) method~\cite{Pol64} is given by
\begin{align*}
	y'_k = x_k + \beta'_k(x_k-x_{k-1})
	\quad \mbox{and} \quad
	x_{k+1} = y'_k-\gamma\grad{f}(x_k),
\end{align*} 
where $\beta_k, \beta'_k >0$ are extrapolation coefficients.  Setting $m_k=1$ in~\eqref{eq:ya}, the AA-GD method is analogous to AGD and HB with $\beta_k, \beta'_k$ replaced by $-\gam{1}{k}$. However, their update directions are chosen differently: AGD takes a step using the gradient at the extrapolated point $y_k$, HB uses the gradient at $x_k$, while AA-GD uses a combination of the gradients evaluated at $x_k$ and $x_{k-1}$. 
For $m>1$, AA-GD is similar to the MiFB method in~\cite{LFP16}. However, unlike AA,  there is currently no efficient way to select the coefficients in~MiFB, thereby restricting its history parameter to $m=1$ or $2$.

\section{Anderson Acceleration for Proximal Gradient Method}\label{sec:pga}

Consider a composite convex minimization problem of the form
\begin{align}\label{eq:comp:optim}
	\underset{x\in \R^n}{\mbox{minimize}} \,\, \varphi(x) := f\left(x\right) + h(x),
\end{align}
where $f: \R^n \to \R$ is $L$-Lipschitz smooth, \emph{i.e.}
\begin{align*}
	\norm{\grad{f(x)}-\grad{f(y)}}_2 \leq L \norm{x-y}_2, \quad \forall x,y \in \dom f.
%	,\\
%\end{align*}
%and
%\begin{align*}
%	f(x) \geq f(y) + \grad{f(y)}^\top(x-y) + \frac{\mu}{2}\norm{x-y}_2^2 , \quad \forall x,y \in \R^n,
\end{align*}
and $h$ is a proper closed and convex function. 
Recall that the proximal operator associated with $h$ is defined as 
\begin{align*}%\label{eq:prox}
	\prox{h}{y} := \argmin_{x} \left\{ h(x) + \frac{1}{2}\norm{x-y}_2^2 \right\}.
\end{align*}
A classical method for solving~\eqref{eq:comp:optim}  is the proximal gradient algorithm (PGA)
\begin{align}\label{eq:PGA:o}
%	y_{k+1} &= x_k-\gamma\grad{f(x_k)}\\\label{eq:PGD:b}
%	x_{k+1} &= \prox{\gamma h}{y_{k+1}}.
	x_{k+1} = \prox{\gamma h}{x_k-\gamma\grad{f(x_k)}},
\end{align}
which can be seen as the fixed-point iteration for the mapping 
\begin{align}\label{eq:g:proxmap}
g(x)= \prox{\gamma h}{x-\gamma\grad{f(x)}}.
\end{align}
 It is not difficult to show that $x^\star$ is a minimizer of~\eqref{eq:comp:optim} if and only if
\begin{align}\label{eq:opt:cond}
	x^\star = \prox{\gamma h}{x^\star-\gamma\grad{f(x^\star)}},
\end{align}
which implies that finding $x^\star$ amounts to finding a fixed-point of $g$. 

In light of our previous discussion, it would be natural to speed-up the PGA method by applying AA to the mapping $g$ in \eqref{eq:g:proxmap}. 
However, in many cases, the function $h$ does not have full domain; for example, when $h$ is the indicator function of some closed convex set. As AA forms an affine (and not a convex) combination in each step, the
resulting iterates can the lie outside $\dom h$ (at which $\nabla f$ may not exist).
Nevertheless, if we  rewrite the PGA iteration as
%on the form:
\begin{align}\label{eq:PGA}
	y_{k+1} = x_k-\gamma\grad{f(x_k)}
	\quad \mbox{and} \quad
	x_{k+1} = \prox{\gamma h}{y_{k+1}},
\end{align}
and consider  the mapping $g$ defined as
\begin{align}\label{eq:g}
	g(y)= \prox{\gamma h}{y}-\gamma\grad{f(\prox{\gamma h}{y})},
\end{align}
then the fixed-point iteration $y_{k+1}=g(y_{k})$ recovers exactly the PGA iteration in \eqref{eq:PGA}. It is clear that  if $y^\star$ is a fixed-point of $g$, then $x^\star=\prox{\gamma h}{y^\star}$ is an optimal solution to~\eqref{eq:comp:optim} since it satisfies condition~\eqref{eq:opt:cond}. 
Now, to relate the convergence of the primal sequence $\{x_k\}$ and the auxiliary $\{y_k\}$, we use the following simple but useful observation: Suppose that $x^\star$ satisfies \eqref{eq:opt:cond}, then $y^\star=x^\star-\gamma\grad{f(x^\star)}$ is a fixed-point of $g$ defined in \eqref{eq:g} and
\begin{align*}
	\norm{x_{k}-x^\star}_2
	=
		\norm{\prox{\gamma h}{y_{k}}-\prox{\gamma h}{x^\star-\gamma\grad{f(x^\star)}}}_2
	\leq
		\norm{y_{k}-y^\star}_2,
\end{align*}
where the last step follows from the nonexpansiveness of proximal operators.
The inequality implies that if one can quickly drive  $\{ y_{k}\}$ to $y^\star$, then $\{x_k\}$ will quickly converge to $x^\star$. It turns out that working with this $g$ is also convenient in designing our safeguarding scheme later.

We thus propose to use AA for accelerating the auxiliary sequence $\{y_k\}$ governed by $g$ defined in~\eqref{eq:g}. Since there are no restrictions on $\{y_k\}$, AA-PGA avoids the feasibility problems \added{of na{\"i}ve AA}. 
Just like PGA, the algorithm requires only one gradient and one proximal evaluation per step. 
The resulting scheme, which we call AA-PGA,  is summarized in Algorithm~\ref{alg:aa:pga}.

\begin{algorithm}[!t]
	\caption{AA-PGA}
	\begin{algorithmic}[1]\label{alg:aa:pga}
		\REQUIRE $x_0=y_0$, $m\geq 0$		
		\STATE $y_1 \gets x_0-\gamma\grad{f(x_0)}$,	$x_1 \gets \prox{\gamma h}{y_1}$, $g_0 \gets y_1$
		\FOR{$k=1,\ldots, K-1$}
		\STATE $m_k \gets \min(m,k)$
		\STATE $g_k \gets x_k - \gamma \grad{f(x_k)}$ and $r_k \gets g_k - y_k$
		\STATE  $R_k \gets [r_k, \ldots, r_{k-m_k}]$ 
		\STATE $ \alpha^k	\gets \argmin_{\alpha^\top\ones=1}\norm{R_k \alpha}$
		\STATE $y_{k+1} \gets  \sum_{i=0}^{m_k} \alpha_i^k g_{k-i}$
		\STATE $x_{k+1} \gets \prox{\gamma h}{y_{k+1}}$
		\ENDFOR
		\ENSURE $x_K$
	\end{algorithmic}
\end{algorithm}	

\subsection{Local Convergence Guarantees}

%Having addressed the feasibility issue in the previous section, 
%our next goal is to derive convergence guarantees for AA-PGA. To this end, we will extend the theory for limited-memory AA to a class of non-smooth problems which includes our proposed algorithm.

Although convergence properties of AA for linear mappings with \emph{full memory} ($m=\infty$) are relatively well understood \cite{WN11, PE13}, much less is known in the case of nonlinear mappings and \emph{limited-memory}. 
The work \cite{TK15} was the first to show that no matter what value $m\in\{0,1,\ldots\}$ is used, AA does not harm the convergence of the fixed-point iteration when started near the fixed point. The proof requires continuous differentiability of $g$. However, in the context of composite convex optimization, the mapping $g$ defined in~\eqref{eq:g} is, in general, non-differentiable. Therefore, the analysis in~\cite{TK15} is not applicable anymore. 
To circumvent this difficulty, we rely on the notion of generalized second-order differentiability, defined below. The interested reader is referred to \cite[Section~13]{RW07} for a comprehensive treatment of epi-differentiability.
\begin{definition}
A function $f$ is twice epi-differentiable at $x$ for a vector $v \in\R^n$ if it is epi-differentiable at $x$ and the second-order quotient functions $\Delta^2_{x,v,t} f$ defined by 
\begin{align*}
	\Delta^2_{x,v,t} f(x') = \left[f(x+ tx') - f(x) - t\InP{v}{x'}\right]/(t^2/2) \quad \mbox{for} \,\,\, t>0,
\end{align*}
epi-converge to a proper function as $t\to 0$. The limit, denoted by $\Delta^2_{x,v} f$,  is then the \replaced{\emph{second-oder epi-derivative of $f$}}{\emph{second epi-derivative} function}.
\end{definition}
 
%Note that when needed, twice epi-differentiability of $h$ will be required only at critical points $x\opt$, and only for the gradient $-\grad{f(x\opt)}$.
We make the following assumption.
\begin{assumption}\label{assumption:local:objective}
Let $x\opt \in \argmin_{x} \varphi(x)$. We assume that:
\begin{enumerate}[label=\rm{(\ref{assumption:local:objective}.}\roman*), leftmargin=*]
\item \label{assumption:pos:hessian}
the function $f$ is of class $\mc{C}^2$ around $x\opt$ and there exists a real $\nu>0$ such that
\begin{align*}
	\nabla^2 f(x\opt) \succeq \nu \IM{}.
\end{align*}

\item  \label{assumption:reg:epidiff}
the convex function $h$ is twice epi-differentiable at $x\opt$ for $-\grad{f(x\opt)}$ and the corresponding second-order epi-derivative is generalized quadratic:
\begin{align*}
	\Delta^2_{x\opt,-\grad{f(x\opt)}} h(\xi)   
	= 
		\begin{cases}
			\frac{1}{2} \InP{\xi}{Q\xi}, & \xi \in L\\
			\infty, & \mbox{otherwise},
		\end{cases} 
\end{align*}
where $L$ is a linear subspace of $\R^n$ and $Q\in\R^{n\times n}$ is a symmetric matrix.
\end{enumerate}
\end{assumption}

Twice epi-differentiable functions, introduced by Rockafellar in~\cite{Roc88}, are remarkable in the sense that they may be both non-smooth and extended real-valued, but still have useful second-order properties. 
%Some remarks on Assumption \ref{assumption:local:objective} are in order. 
%Twice epi-differentiable functions,  introduced by Rockafellar in \cite{Roc88}, are remarkable for their second-order properties in variational analysis 
%introduced by Rockafellar in \cite{Roc88} and further developed in \cite{PR96, PR95}. This property is
%and \replaced{include}{enjoyed by} all 
One important class of twice epi-differentiable functions are known as \emph{fully amenable}~\cite{PR95}.
In the context of (additive) composite optimization, full amenability is justified whenever $f\in\mc{C}^2$ and  $h$ is a \emph{polyhedral} function (i.e., its epigraph is a polyhedral set). Indeed, \cite[Proposition~2.6]{PR95} ensures that $\varphi = f + h$ is fully amenable at any feasible $x$, which in turn implies twice epi-differentiability of $h$ at $x$ for $-\grad{f(x)}$ since $\partial \varphi(x) = \grad{f(x)} + \partial h(x)$.
Notable examples of polyhedral $h$ in machine learning applications are the $\ell_1$-norm, $\ell_\infty$-norm, total variation seminorm, and the indicator functions of  polyhedral sets such as the non-negative orthant, box constraints and the probability simplex. 
%Secondly, the case in which $\Delta^2_{x\opt,-\grad{f(x\opt)}} h$ is generalized quadratic is rich enough to... 

For the preceding $\varphi$, it is shown in \cite[Proposition~4.12]{PR95},  that the function $\Delta^2_{x\opt,-\grad{f(x\opt)}} h$ is generalized quadratic if and only if $(x\opt, \grad{f(x\opt)})$ satisfies the \emph{non-degeneracy} condition:
\begin{align}\label{assumption:cq}
	-\grad{f(x\opt)} \in \relint(\partial h(x\opt)).
\end{align}
More broadly, if condition~\eqref{assumption:cq} holds, then any ${\mathcal C}^2$-partly smooth function $h$ satisfies the properties in~\ref{assumption:reg:epidiff} (this follows by combining~\cite[Theorem~28]{DHM06} and \cite[Theorem~4.1(a) and (g)]{PR96};  see \cite{TSP18} for detailed  arguments). This allows to  include regularizers which are not polyhedral, like the nuclear norm in matrix completion and the $\ell_1\!\!-\!\ell_2$-norm in group lasso \cite{LFP17}. 

%Thirdly, we remark that  the class of twice epi-differentiable functions has been significantly enlarged recently, which covers constraint sets, which go far beyond polyhedral ones,  such as the second order cone and more \cite{HMS20, MMS19}.

Note that condition \eqref{assumption:cq} is very mild and  can be seen as a geometric generalization of the well-known \emph{strict complementarity} in nonlinear programming \cite{Bur90}. For example, for the lasso problem with $f(x)= (1/2)\ltwo{Ax-b}^2$ and $h(x) = \lambda \lone{x}$, it is easy to verify that \eqref{assumption:cq} is justified as long as $|\left(A^\top (Ax\opt - b)\right)_i|\neq \lambda$ whenever $(x\opt)_i=0$. In fact,  this condition has been considered almost necessary for identifying the support of $x\opt$ \cite{LFP17}. 

An important consequence of Assumption \ref{assumption:local:objective} is that the proximal mapping $\proxmap_{\gamma h}$ becomes differentiable at $y\opt=x\opt-\gamma \grad{f(x\opt)}$. This fact is summarized in the following lemma.
\begin{lemma}
	Let Assumption \ref{assumption:local:objective} hold. Then, the proximal operator $\proxmap_{\gamma h}$ is  differentiable at $y\opt=x\opt - \gamma \grad{f(x\opt)}$ and its  Jacobian $P_{\gamma}(y\opt) := \mathop{J}_{\proxmap_{\gamma h}} \left(y\opt\right)$ is symmetric and positive semidefinite with $\ltwo{P_{\gamma}(y\opt)} \leq 1$. Moreover,  
	the mapping $g$ is  differentiable at $y\opt$ with Jacobian:
	\begin{align*}
	G = P_{\gamma}(y\opt)\left(\IM{}  - \gamma   \nabla^2 f(x\opt)\right).
	\end{align*}
	If, in addition, $\gamma\in ( 0, 1/L]$, then  $\ltwo{G}  \leq 1 - \gamma \nu \in [0,1)$. 
\iffalse
\begin{lemma}
	Let Assumption \ref{assumption:local:objective} hold. Then, the proximal operator $\proxmap_{\gamma h}$ is  differentiable at $y\opt=x\opt - \gamma \grad{f(x\opt)}$ with  Jacobian $P_{\gamma}(y\opt) := \mathop{J}\proxmap_{\gamma h} \left(y\opt\right)$, 
	and the mapping $g$ is  differentiable at $y\opt$ with Jacobian:
	\begin{align*}
	G = P_{\gamma}(y\opt)\left(\IM{}  - \gamma   \nabla^2 f(x\opt)\right).
	\end{align*}
	Moreover, $P_{\gamma}(y\opt)$ is symmetric and positive semidefinite with $\ltwo{P_{\gamma}(y\opt)} \leq 1$. If furthermore $\gamma\in ( 0, 1/L]$, we also have $\ltwo{G} \leq 1 - \gamma \nu\in [0,1)$. 
	\begin{proof}
\fi
	\begin{proof}
Detailed arguments for differentiability of $\proxmap_{\gamma h}$ at $y\opt$ can be found in \cite[Thm.~4.10]{TSP18}. The Jacobian $G$ of $g$ at $y\opt$ is a direct consequence of the chain rule. Finally, since $\nabla^2 f(x\opt) \succeq \nu \IM{}$ and $f$ is $L$-smooth, we have $\ltwo{G}\leq \ltwo{P_{\gamma}(y\opt)} \ltwo{\left(\IM{}  - \gamma   \nabla^2 f(x\opt)\right)} \leq 1-\gamma\nu$, as desired.
\end{proof}
\end{lemma}

Our last assumption imposes a boundedness condition on the extrapolation coefficients.
\begin{assumption}\label{assumption:bounded:extra:coeffs}
There exists a constant $M_\alpha$ such that $\norm{\alpha^k}_1 \leq M_\alpha$ for all $k\in \N_+$. 
\end{assumption}
This assumption is very common in the literature of AA and some effective solutions have been proposed to enfore it in practice. For example, one can monitor the condition number of the $R$ matrix in the QR decomposition and drop the left-most column of the matrix if the number becomes too large \cite{WN11}, or one can add a Tikhonov regularization to the least squares as was done in \cite{SdB16}. 
The condition can also be imposed directly in the algorithm without changing the subsequent results. More specifically, if \added{we detect that} $\norm{\alpha^k}_1$ is greater than $M_\alpha$, \replaced{we can}{one can}  set $\alpha^k=\left[0,.\ldots,1\right]^\top$, i.e., we simply perform a fixed-point iteration step.

We can now state the main result of this section.
\begin{theorem}\label{thrm:main}
Let Assumptions  \ref{assumption:local:objective} and \ref{assumption:bounded:extra:coeffs} hold. Let $\gamma \in (0, 1/L]$ and define $\rho(G)=\ltwo{G}$. Let  $\hat{\rho}$ be some real constant satisfying $\hat{\rho} \in (\rho(G),1)$. Let $F(y) = g(y)-y$ with $g$ given in \eqref{eq:g} and let $y^\star = x\opt - \gamma \grad{f(x\opt)}$ be a fixed-point of $g$. If $y_0$ is initialized sufficiently close to $y^\star$, then, for any fixed $m\in \N$, the iterates $\{x_k\}$ and $\{y_k\}$ formed by AA-PGA satisfy:
\begin{align*}
		\norm{F(y_k)}_2 \leq \hat{\rho}^k\norm{F(y_0)}_2
%		\\
	\quad \mbox{and} \quad
	 \norm{x_k-x^\star}_2 \leq \frac{3+\rho(G)}{1-\rho(G)}\hat{\rho}^k\norm{y_0-y^{\star}}_2.
\end{align*}
Moreover, we have
\begin{align*}%\label{eq:gy:main:thrm}
	\limsup_{k\to \infty} \left(\frac{\ltwo{F(y_k)}}{\ltwo{F(y_0)}}\right)^{1/k} \leq \rho(G)
	\quad\mbox{and}\quad
	\limsup_{k\to \infty} \left(\frac{\ltwo{x_k-x\opt}}{\ltwo{x_0-x\opt}}\right)^{1/k} \leq \rho(G).
\end{align*}
\begin{proof}
See Appendix~\ref{appendix:thrm:main}.
\end{proof}
\end{theorem}

The theorem implies that when initialized near the optimal solution, even in the worst case, the use of multiple past iterates to construct a new update in AA will not slow down the convergence of the original PGA method, no matter how we choose $m\in\{0,1,\ldots\}$. In most cases, near the solution, we would expect AA-PGA to enjoy the strong adaptive rate in \eqref{eq:opt:rate} even for a small value of $m$. Therefore, we can see AA as interpolating between the two convergence rates corresponding to $m=0$ (PGA) and $m=\infty$ (full-memory AA). Whether or not AA can attain a stronger convergence rate guarantees than PGA for finite $m$ is still an open question, even with smooth and linear mappings. 

\section{Guarded Anderson Accelerated PGA}

We have shown that when started from a point close to the  optimal solution, AA-PGA is convergent under mild conditions. 
A natural question, which has also recently been raised in \cite{FZB19}, is whether AA converges globally. We show that the answer is negative even when the problem has no constraint and the objective function is smooth. In this case, AA-PGA reduces to AA-GD, and hence the result is also valid for the AA methods in \cite{WN11,SdB16}.  To that end, we construct a one-dimensional smooth and strongly convex function and show analytically that AA will not converge to the optimum but get stuck in a periodic orbit. Concretely, consider the function $f$ whose gradient is given by
\begin{align}\label{eq:aa:myfunc}
	\grad{f(x)} = \begin{cases}
		\frac{x}{10} - 24.9 \quad & \mbox{if} \quad x<-1,\\
		25x \quad & \mbox{if} \quad -1 \leq x < 1,\\
		\frac{x}{10} + 24.9 \quad & \mbox{if} \quad x\geq 1.
	\end{cases}
\end{align}
This $f$ is strongly convex with $\mu=1/10$ and smooth with $L=25$. 
\begin{figure}[h!]
	\centering 
	\begin{minipage}{0.5\textwidth}
		\centering 
		{\includegraphics[width=\textwidth]{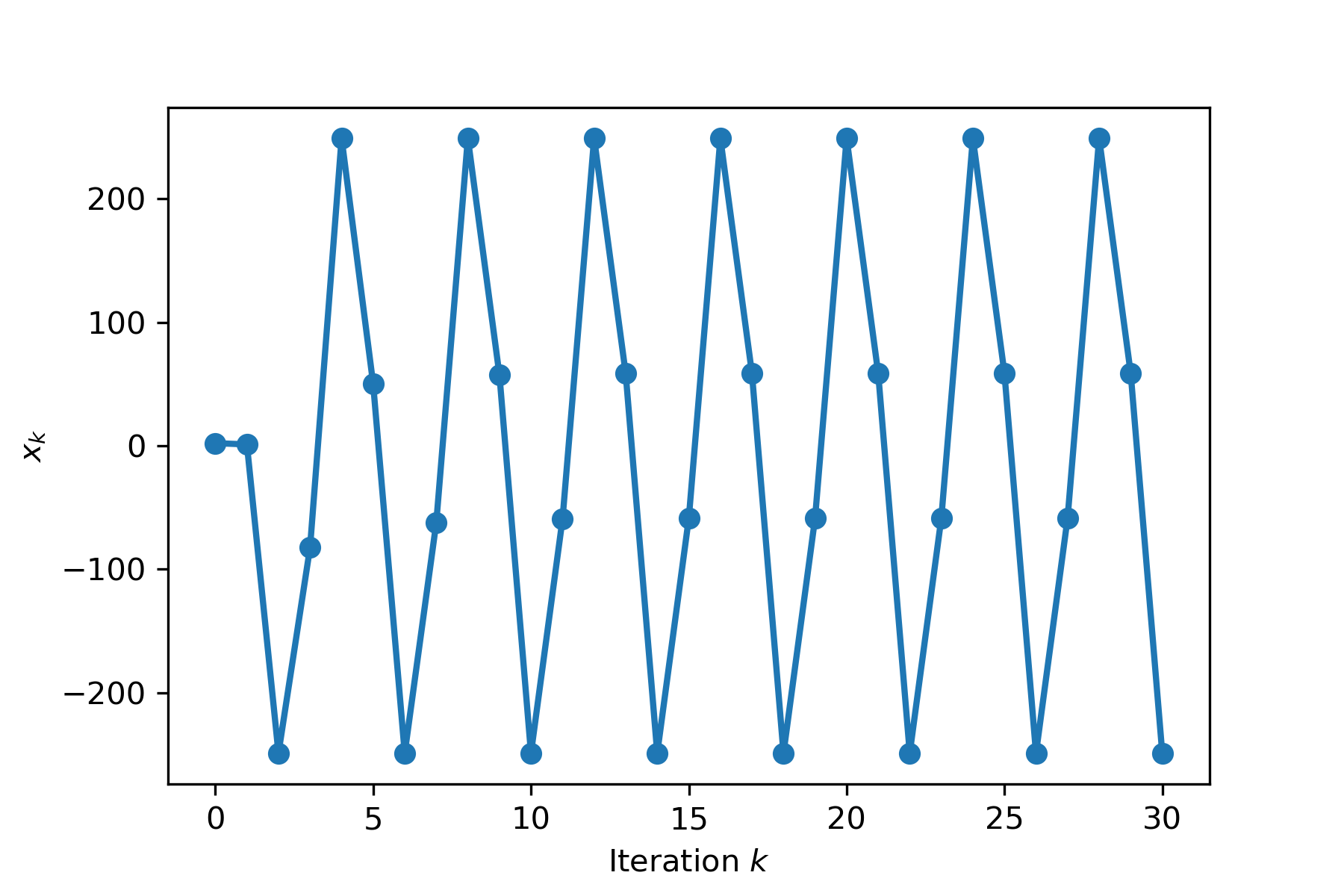}}	
	\end{minipage}
	\hskip -0.15in
	\begin{minipage}{0.5\textwidth}
		\centering 
	{\includegraphics[width=\textwidth]{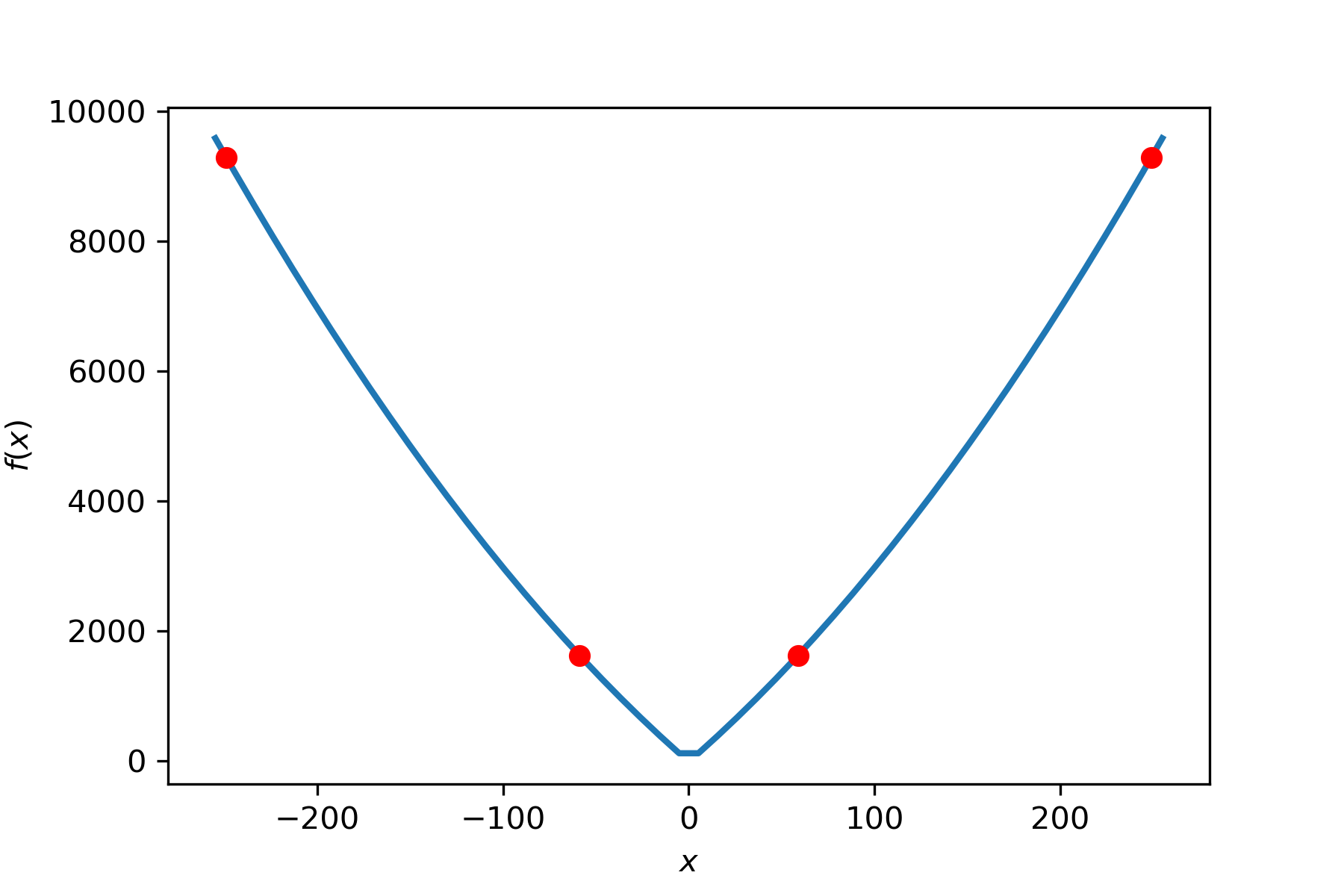}}	
	\end{minipage}
	\caption{Left: Iterates of the AA-GD method when minimizing $f(x)$ defined in \eqref{eq:aa:myfunc} with $x_0=2.1$. Right: The graph of $f(x)$ with the circles indicating the four limit points shown in Proposition~\ref{prop:example}.} \label{fig:aa:limit:cycle}
%	\vskip -0.075in
\end{figure}
%We consider AA-GD with $m=1$ for which the iterates can be expressed explicitly as (see Appendix~\ref{appendix:prop1}):
%\begin{align}\label{eq:aa:example:catsys}
%	\begin{bmatrix}
%	x_{k+1}\\
%	x_{k}
%	\end{bmatrix}
%	=
%	\begin{bmatrix}
%		\frac{\grad{f(x_{k-1})}}{\grad{f(x_{k-1})}- \grad{f(x_k)}} & \frac{-\grad{f(x_k)}}{\grad{f(x_{k-1})}- \grad{f(x_k)}}\\
%		1 & 0
%	\end{bmatrix}
%	\begin{bmatrix}
%		x_{k}\\
%		x_{k-1}
%		\end{bmatrix}.
%\end{align}
%\begin{align*}
%	x_{k+1} 
%	= 
%		\frac{\grad{f(x_{k-1})}x_k }{\grad{f(x_{k-1})}- \grad{f(x_k)}}
%		-
%		\frac{\grad{f(x_k)}x_{k-1}}{\grad{f(x_{k-1})}- \grad{f(x_k)}}.
%\end{align*}
A trajectory of AA-GD with $m=1$ started at $x_0=2.1$ is depicted in Figure~\ref{fig:aa:limit:cycle} indicating that it converges to a periodic orbit instead of the origin. More formally, one can show the following.
\begin{proposition}\label{prop:example}
Let $f$ be the function defined in \eqref{eq:aa:myfunc}. Suppose that the AA-GD method is applied to minimize $f$ with the history parameter $m=1$ and the step size $\gamma=1/L$.  Then, for any initial point $x_0\in [2.01, 246.98]$ and $n = 0,1,\ldots$, the iterates generated by AA-GD satisfy:
\begin{align*}
	x_{4n+3} \to -249(\sqrt{5}-2), \quad  x_{4n+4} = +249, 
	\\ 
	x_{4n+5} \to +249(\sqrt{5}-2),\quad  x_{4n+6} = -249.
\end{align*}
\begin{proof}
See Appendix~\ref{appendix:prop1}.
\end{proof}
\end{proposition}

The proposition confirms the necessary of a safeguarding step to ensure global convergence (see, e.g., \cite{ZOB18,FZB19}). Note that such a step often only involves checking a simple condition, and hence  is cheaper to execute than a line search.

Recall that each iteration of AA-PGA consists of one original PGA step followed by an AA step. Thus, one natural strategy for stabilization would be to compare the objective value produced by the AA step with that of PGA and select the one with lower value as the next iterate. 
%Since PGA converges globally, this ensures that AA-PGA also converges globally at least as fast as PGA.  
However, this approach can be costly since one needs two function evaluations per step. 
Indeed, only the descent condition below is needed to achieve the same convergence rate as PGA: 
\begin{align}\label{eq:desc:cond}
	f(x_{k+1}) 
%	\leq 
%		f(x_k) + \InP{\grad{f(x_k)}}{x_{k+1}-x_k} + \frac{1}{2\gamma}\norm{x_{k+1}-x_k}_2^2
	\leq
	 	f(x_k) - \frac{\gamma}{2}\norm{\grad{f(x_k)}}_2^2.
\end{align}

\begin{algorithm}[h]
	\caption{Guared AA-PGA}
	\begin{algorithmic}[1]\label{alg:aa:pga:gd}
		\REQUIRE $y_0 = x_0$, $m\geq 0$		
		\STATE $y_1 \gets x_0-\gamma\grad{f(x_0)}$,	$x_1\gets \prox{\gamma h}{y_1}$, $g_0 \gets y_1$, and $r_0 = g_0 - y_0$
		\FOR{$k=1,\ldots, K-1$}
			\STATE 	$m_k \gets \min(m,k)$
			\STATE $g_k \gets x_k - \gamma \grad{f(x_k)}$ and $r_k \gets g_k - y_k$
			\STATE  $R_k \gets [r_k, \ldots, r_{k-m_k}]$
			\STATE $ \alpha^k	\gets \argmin_{\alpha^\top\ones=1}\norm{R_k \alpha}$
			\STATE $y_{\mathrm{ext}} \gets  \sum_{i=0}^{m_k} \alpha_i^k g_{k-i}$
			\STATE $x_{\mathrm{test}} \gets \prox{\gamma h}{y_{\mathrm{ext}}}$
			\IF{$f(x_{\mathrm{test}})  \leq  f(x_k) -\frac{\gamma}{2}\norm{\grad{f(x_k)}}_2^2$} 
				\STATE $x_{k+1}=x_{\mathrm{test}}$ \,\, and \,\, $y_{k+1}=y_{\mathrm{ext}}$ 
%				\STATE $f(x_{k+1})=f(x_{\mathrm{test}})$
			\ELSE 
				\STATE $x_{k+1}=\prox{\gamma h}{g_k}$ \,\, and \,\, $y_{k+1}=g_k$ 						
			\ENDIF
		\ENDFOR
		\ENSURE $x_K$
	\end{algorithmic}
\end{algorithm}	

This suggests an alternative way for stabilization, which is to compare the objective value of the AA step with the right-hand side of \eqref{eq:desc:cond}. If sufficient descent was made, the AA step is accepted, otherwise the PGA step is chosen. This allows to reuse the function values more efficiently. In particular, if the AA step is selected, only one function evaluation is needed. Moreover, in many applications, function values can be computed at a very small additional cost by reusing information readily available from gradient evaluations. 
Putting everything together, we arrive at Algorithm~\ref{alg:aa:pga:gd} that admits the global convergence rate of PGA with the potential for local adaptation and acceleration.

\begin{proposition}[Global convergence]
Let $f:\mathbb{R}^n\mapsto\mathbb{R}$ be  $\mu$-strongly convex and $L$-smooth  and let $\gamma\in(0,2/(\mu+L)]$. Then, the iterates $\{x_k\}$  generated by Algorithm~\ref{alg:aa:pga:gd} satisfy
\begin{align}
	\norm{x_k-x^\star}_2^2 \leq O\left(1-\frac{\gamma\mu L}{\mu+L}\right)^k\norm{x_0-x^\star}_2^2.
\end{align}
\end{proposition}
The proof of this result is straightforward (see, e.g., \cite{Nes04, Bec17}), and  follows directly from the descent condition \eqref{eq:desc:cond} and a standard strong convexity inequality. Hence, we omit it here.

\section{Extension to Bregman proximal gradient methods}
Consider optimization problems of the form
\begin{align}\label{eq:breg:optim}
	\underset{x\in \mathcal{D}}{\mbox{minimize}} \,\, f\left(x\right) + h(x),
\end{align}
where $\mathcal{D} \subseteq \R^n$ is a closed convex set with nonempty interior. The formulation~\eqref{eq:breg:optim} often provides a more flexible way to handle the constraints, which are usually encoded by $h$ in \eqref{eq:comp:optim}. This model is very rich and led to several recent advances in algorithmic developments of first-order methods. Bregman proximal gradient (BPG) is a general and powerful tool for solving \eqref{eq:breg:optim} thanks to its ability to exploit the underlying geometry of the problem.  
The mirror descent method \cite{NY83,BT03} is a well-known instance of BPG when $h(x)=\mathrm{I}_{\mathcal{C}}(x)$ for some closed convex set $\mathcal{C} \subseteq \mathcal{D}$. Some more recent instances of BPG include the NoLips algorithm \cite{BBT16} and its accelerated version analysed in \cite{HRX18}. The number of applications of the BPG framework are growing rapidly~\cite{BSTV18, DBd19,LFN18}.

The BPG method fits the geometry of the problem at hand, which is typically governed by the constraints and/or the objective,  all-in-one by means of a  kernel function. Popular examples include the energy function $\genf(x)=(1/2)\norm{x}_2^2$; the Shannon entropy $\genf(x) = \sum_{i=1}^{n}x_i \log x_i$, $\dom \genf = \R^n_+$ with ($0\log 0=0$); the Burg entropy $\genf(x)=-\sum_{i=1}^{n} \log x_i$, $\dom \genf = \R^n_{++}$; the Fermi-Dirac entropy $\genf(x)=\sum_{i=1}^{n}\left( x_i\log x_i + (1-x_i)\log (1-x_i)\right)$, $\dom\genf=[0,1]^n$; the Hellinger entropy $\genf(x)= -\sum_{i=1}^{n} \sqrt{1-x_i^2}$, $\dom \genf = [-1,1]^n$; and the polynomial function $\genf(x)=\frac{\alpha}{2}\norm{x}_2^2+\frac{1}{4}\norm{x}_2^4$, $\alpha\geq 0$.

We impose the following assumption in this section.
%\begin{assumption}
%The function $\genf: \R^n \to (-\infty, +\infty]$ is of Legendre type with $\overline{\dom} \genf=\mathcal{D}$ and its conjugate $\genf^*$ satisfies $\dom \grad{\genf^*}=\R^n$.
%\end{assumption}
%We make the following additional assumptions in this section.
%\begin{assumption}
% The set $\overline{\dom} \genf=\mathcal{D}$ is convex and the following conditions hold:
%\begin{description}
%	\item [1.]  $\genf: \R^n \to (-\infty, +\infty]$ is of Legendre type  and its conjugate $\genf^*$ satisfies $\dom \grad{\genf^*}=\R^n$.
%	\item [1.] $f: \R^n \to (-\infty, +\infty]$ is proper closed convex  and differentiable on $\intdom \genf$.	
%	\item [2.] $h: \R^n \to (-\infty, +\infty]$ is proper closed convex and $\dom h \cap \intdom \genf \neq \emptyset$.
%\end{description}
%\end{assumption}
\begin{assumption}\label{assumption:2}
 The set $\overline{\dom} \genf=\mathcal{D}$ is convex and the following conditions hold:\\[1.ex]
1.  $\genf: \R^n \to (-\infty, +\infty]$ is of Legendre type  and its conjugate $\genf^*$ satisfies $\dom \grad{\genf^*}=\R^n$.\\[0.5ex]
2.	 $f: \R^n \to (-\infty, +\infty]$ is proper closed convex  and differentiable on $\intdom \genf$.	\\[0.5ex]
3.  $h: \R^n \to (-\infty, +\infty]$ is proper closed convex and $\dom h \cap \intdom \genf \neq \emptyset$.
\end{assumption}

When $\genf$ is Legendre, its gradient $\grad{\genf}$ is a bijection from $\intdom \genf$ to $\intdom \genf^*$ while $\grad{\genf^*}$ is a bijection from $\intdom \genf^*$ to $\intdom \genf$, i.e., $(\grad{\genf})^{-1}=\grad{\genf^*}$ \cite[Chapter~26]{Roc70}. 
%Note also that ${\rm ran} \grad{\genf}=\dom \grad{\genf^*}=\intdom \genf^*$. 
Note that in all the above examples, $\genf$ is Legendre. Moreover, except from the Burg entropy, all the others share the useful property $\dom \grad{\genf^*}=\R^n$, which is critical for the development of our AA scheme.

The Bregman distance associated with $\genf$ is the function $D_\genf: \dom \genf \times \intdom \genf \to \R$ given by
\begin{align*}
	\breg{x,y} = \genf(x) - \genf(y) - \InP{\grad{\genf(y)}}{x-y}.	
\end{align*}
At the core of the BPG method is the Bregman proximal operator that generalizes the conventional one and is defined as~\cite{CZ92}:
\begin{align}\label{eq:breg:prox}
	\proxmap_{h}^\genf(y) = \argmin_{x \in \R^n} \left\{h(x) + \breg{x, y}\right\}, \quad y \in \intdom \genf.
\end{align}
BPG starts with some $x_0 \in \intdom \genf$ and performs the following operator at each iteration:
\begin{align*}
	x_{k+1}=T_\gamma(x_k) 
	:= 
		\argmin_{x \in \R^n}
		\big\{
				f(x_k)
				+
				\InP{\grad{f(x_k)}}{x-x_k}+\gamma^{-1} \breg{x, x_k} 
				+
				h(x)
		\big\}.
\end{align*}
Assumptions~\ref{assumption:2} ensures that BPG iterates are well-defined and $x_k \in \intdom \genf$ for all $k$~\cite[Lemma~2]{BBT16}. 
Further simplification the update formula yields 
\begin{align*}
	x_{k+1} 
	= 
		\argmin_{x \in \R^n} 
			\left\{
				\InP{\gamma\grad{f(x_k)}-\grad{\genf(x_k)}}{x}
				+
				\genf(x)
				+ h(x)
			\right\}.
\end{align*}
Using the optimality condition and the fact that $\left(\grad{\genf}\right)^{-1}=\grad{\genf^*}$ yield
\begin{align}\label{eq:opt:cond:breg:equi}
	0 \in \gamma \partial h(x_{k+1}) + \grad{\genf(x_{k+1})} - \grad{\genf(\grad{\genf^*} \left(\grad{\genf(x_k)} - \gamma \grad{f(x_k)}\right))}.
\end{align}
Comparing \eqref{eq:opt:cond:breg:equi} with the optimality condition of \eqref{eq:breg:prox}, we obtain an equivalent update rule for BPG: 
\begin{align*}
	x_{k+1} = \proxmap_{\gamma h}^\genf\left(\grad{\genf^*} \left(\grad{\genf(x_k)} - \gamma \grad{f(x_k)}\right)\right).
\end{align*}
Note that when $\genf$ is the energy function, $\grad{\genf}$ and $\grad{\genf^*}$ are the identity map\replaced{ and}{,} we recover the PGA method.
To apply AA, we further express the BPG iterations on the form
\begin{align}\label{eq:md:2}
	y_{k+1} = \grad{\genf(x_k)} - \gamma \grad{f(x_k)}
	\quad \mbox{and} \quad
	x_{k+1} = \proxmap_{\gamma h}^\genf\left(\grad{\genf^*} \left(y_{k+1}\right)\right).	
\end{align}
In words, the mirror map $\grad{\genf}$ maps $x_k$ from the primal space to \replaced{a}{the} dual one, where the gradients live. A gradient step is then taken in the dual space to obtain $y_{k+1}$. Next, $y_{k+1}$ is transferred back to the primal space by the inverse map $\grad{\genf^*}$.  Finally, the Bregman proximal operator is performed in the primal space to produce $x_{k+1}$.

Our strategy is to extrapolate the sequence $\{y_k\}$. Note that this sequence can be seen as  the fixed-point iteration of 
\begin{align*}
	g(y) = 
		\grad{
			\genf(
				\proxmap_{\gamma h}^\genf \circ \grad{\genf^*}(y)
			)} 
		-
		 \gamma \grad{f(\proxmap_{\gamma h}^\genf \circ\grad{\genf^*}(y))}.
\end{align*}
The AA scheme applied to this mapping (called AA-BPG) has a simple and elegant interpretation. Concretely, instead of accelerating the primal sequence, which is restricted to the constraint set, it extrapolates a sequence in the dual space, avoiding feasibility issues since $\grad{\genf^*}$ has full domain.
To gain some intuition, we first recall the following useful property of Legendre functions:
\begin{align*}
 \breg{\grad{\genf^*(y)}, \grad{\genf^*(y')}} = D_{\genf^*}\left(y',y\right) \quad \forall y, y' \in \intdom \genf^*.
\end{align*}
Assume that $g$ has a fixed-point $y^\star$ and $\{y_k\}$ generated by AA-BPG is converging to $y^\star$. Let $\grad{\genf^*(y_k)}$  and $\grad{\genf^*(y^\star)}$ be the images of $y_k$ and $y^\star$ on the primal space, then it holds that
\begin{align*}
 D_{\genf^*}\left(y^\star,y_k\right) = \breg{\grad{\genf^*(y_k)}, \grad{\genf^*(y^\star)}}.
\end{align*}
Applying the Bregman operator to the two images will give us $x_{k}$ and $x^\star$, respectively. Since Bregman proximal operators possess certain nonexpansiveness property akin to their Euclidean counterpart~\cite{BI12,Eck93}, it is thus reasonable to expect that $\breg{x_k,x^\star}$ is well approximated by $D_{\genf^*}\left(y^\star,y_k\right)$; for example, when $\dom h \subseteq \intdom \genf$, it is shown in \cite{BI12} that $ \breg{x_k,x^\star} \leq \breg{\grad{\genf^*(y_k)}, \grad{\genf^*(y^\star)}}$. 
Moreover, $y_k = y^\star$ implies $x_k=x^\star$. Therefore, if AA can speed-up the convergence of $\{y_k\}$, one can achieve similar acceleration for $\{x_k\}$. 

In the above discussion, we implicitly assumed that $x^\star\in\intdom \genf$. However, if $x^\star$ happens to be on the boundary of $\dom \genf$, the mirror map $\grad{\genf}$ at $x^\star$ does not exist. One can then no longer express $x^\star$ as a fixed-point of some mapping involving $\grad{\genf}$. This makes it very hard to derive general theoretical guarantees for BPG since essentially all the current proofs of AA are  heavily based on $g(x^\star)$. Therefore, a new proof technique that goes beyond linearization of $g$ around $x^\star$ is needed, which we leave as a topic for future research.
Nonetheless, since each iteration of AA-BPG consists of one BPG step, $T_\gamma(x_k)$, one can always compare the progress made by the AA step with the BPG one as was done in AA-PGA.
A counterpart of the sufficient descent condition~\eqref{eq:desc:cond} that ensures the global convergence  of BPG is~\cite{BBT16}:
\begin{align*}
	f(x_{k+1}) 
	\leq 
		f(x_k) + \InP{\grad{f(x_k)}}{T_\gamma(x_k) -x_k} + \gamma^{-1}\breg{T_\gamma(x_k) , x_k}.
\end{align*}
Thus, a similar policy for stabilization as in AA-PGA will retains the convergence rate of BPG. The final AA-BPG algorithm is reported in Algorithm~\ref{alg:aa:bpg:gd}.
\begin{algorithm}[h]
	\caption{Guared AA-BPG}
	\begin{algorithmic}[1]\label{alg:aa:bpg:gd}
		\REQUIRE $y_0 \in \intdom \genf^*$, $x_0 = \proxmap_{\gamma h}^\genf\left(\grad{\genf^*} \left(y_0\right)\right)$, $m\geq 0$		
		\STATE $y_1 \gets \grad{\genf(x_0)} - \gamma \grad{f(x_0)}$,	
				$x_1\gets x_0 = \proxmap_{\gamma h}^\genf\left(\grad{\genf^*} \left(y_0\right)\right)$, $g_0 \gets y_1$, and $r_0 = g_0 - y_0$
		\FOR{$k=1,\ldots, K-1$}
			\STATE 	$m_k \gets \min(m,k)$
			\STATE $g_k \gets \grad{\genf(x_k)} - \gamma \grad{f(x_k)}$ and $r_k \gets g_k - y_k$
			\STATE  $R_k \gets [r_k, \ldots, r_{k-m_k}]$
			\STATE $ \alpha^k	\gets \argmin_{\alpha^\top\ones=1}\norm{R_k \alpha}$
			\STATE $y_{\mathrm{ext}} \gets  \sum_{i=0}^{m_k} \alpha_i^k g_{k-i}$
			and $x_{\mathrm{test}} \gets \proxmap_{\gamma h}^\genf\left(\grad{\genf^*} \left(y_{\mathrm{ext}}\right)\right)$
			\STATE $x_{\mathrm{BPG}} = \proxmap_{\gamma h}^\genf\left(\grad{\genf^*} \left(g_k\right)\right)$
			\IF{$f(x_{\mathrm{test}})  \leq  f(x_k) + \InP{\grad{f(x_k)}}{x_{\mathrm{BPG}} -x_k} + \gamma^{-1}\breg{x_{\mathrm{BPG}} , x_k}$} 
				\STATE $x_{k+1}=x_{\mathrm{test}}$
				\quad and \quad $y_{k+1}=y_{\mathrm{ext}}$ 
%				\STATE $f(x_{k+1})=f(x_{\mathrm{test}})$
			\ELSE 
				\STATE $x_{k+1}=x_{\mathrm{BPG}}$
				\quad and \quad $y_{k+1}=g_k$ 						
			\ENDIF
		\ENDFOR
		\ENSURE $x_K$
	\end{algorithmic}
\end{algorithm}	

\section{Numerical Experiments}
We will now illustrate the performance of (guarded) AA-PGA and AA-BPG on several constrained optimization problems with important applications in signal processing and machine learning. All the experiments are implemented in Python and run on a laptop with four 2.4 GHz cores and 16 GB  of RAM, running Ubuntu 16.04 LTS.

For AA-PGA, we compare it with PGA, PGA with adaptive line search (PGA-LS), and accelerated PGA (APGA) \cite{Bec17}.  For AA-BPG, we compare AA-BPG with BPG, accelerated BPG (ABPG), ABPG with adaptive line search (ABPG-g), and restarted ABPG (ABPG-Restart)~\cite{HRX18}. 
%For all the methods, except from the ones using line search, we set $\gamma=1/L$, where $L=\norm{A}_2^2/4m$ and $L=\norm{A}_2^2/m$,  for the first and the second problems, respectively.
For the AA schemes, we use  $m=5$ in all plots and simply add a Tikhonov regularization of $10^{-10}\norm{R_k}_2^2$ to~\eqref{eq:ls} to avoid singularity, as was done in~\cite{SdB17}, without any tunning. For each experiment, we plot the errors, defined as $f(x_k)-f(x^\star)$, versus the number of iterations and wall-clock runtime. We have picked a few real-world data sets, which are known to be very ill-conditioned, and  hence  challenging for any first order methods.\footnote{The data sets Madelon and Gisette  are downloaded from: \url{http://archive.ics.uci.edu/ml/datasets}. The data sets Cina0 and Sido0  are downloaded from: \url{http://www.causality.inf.ethz.ch}}
 All methods are initialized at $x_0=\textbf{0}$ unless otherwise stated. 

\subsection{Constrained logistic regression} 
We start our experiments with the logistic regression with bounded constraint:
\begin{align*}
	&\underset{x\in \R^n}{\mbox{minimize}} \,\, \frac{1}{M} \sum_{i=1}^{M}\log(1+ \exp(-y_i a_i^\top x)) + \mu \norm{x}_2^2\\
	&\mbox{subject to}\,\, \norm{x}_\infty \leq 1,
\end{align*}
where $a_i\in\R^n$ are training samples  and $y_i \in\{-1,1\}$ are the corresponding labels. We set $\gamma=1/L$, where $L=\norm{A}_2^2/4M$ with $A=[a_1, \ldots, a_M]$.

\begin{figure}[t!]
	\centering 
		\subfigure[Madelon: $\mu=10$, $\kappa=3\times 10^{6}$]{
			{\includegraphics[width=0.4\textwidth]{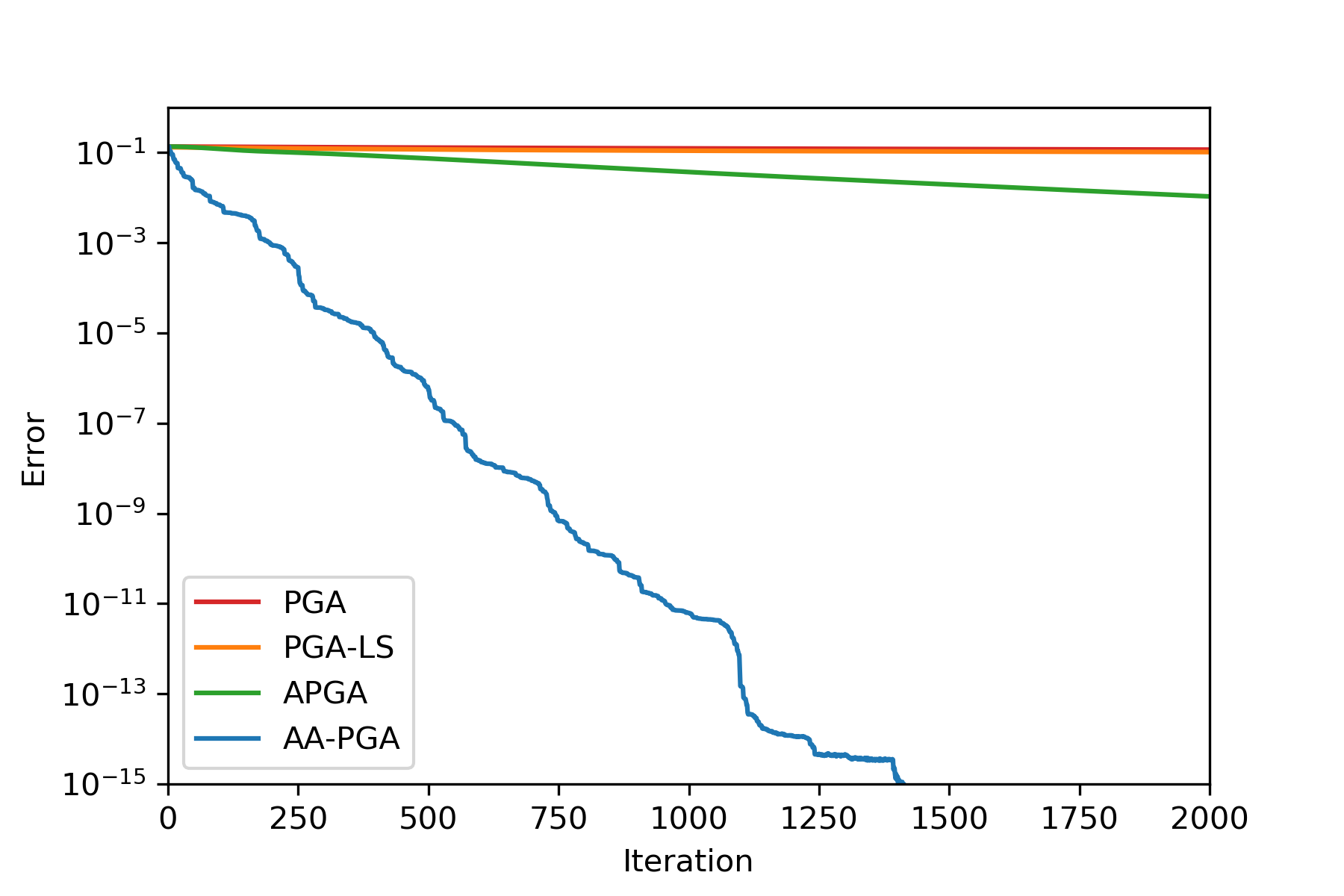}}
			{\includegraphics[width=0.4\textwidth]{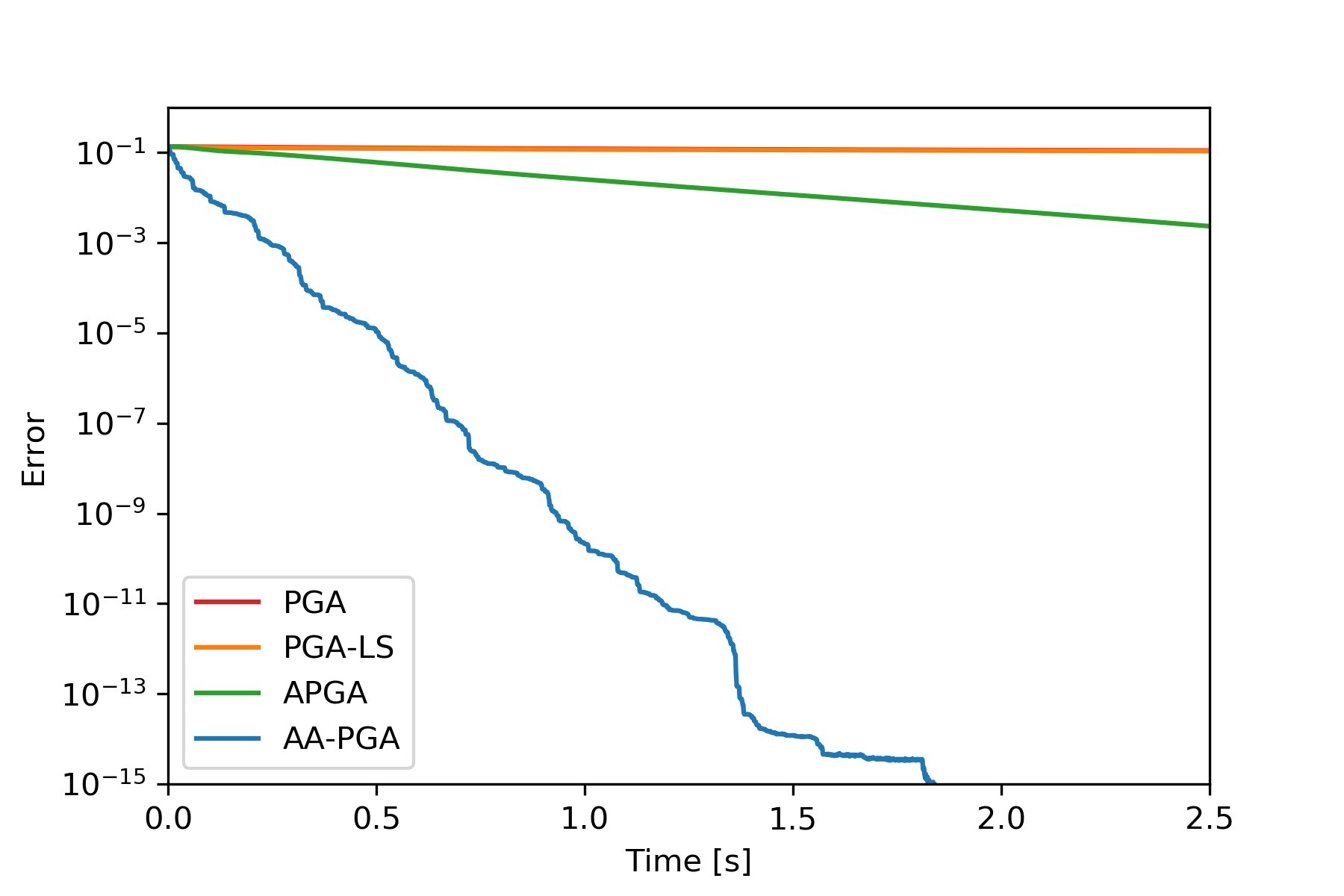}}
		}
	\hfill		
	\subfigure[Gisette: $\mu=10$, $\kappa=3.4\times 10^{6}$]{
			{\includegraphics[width=0.4\textwidth]{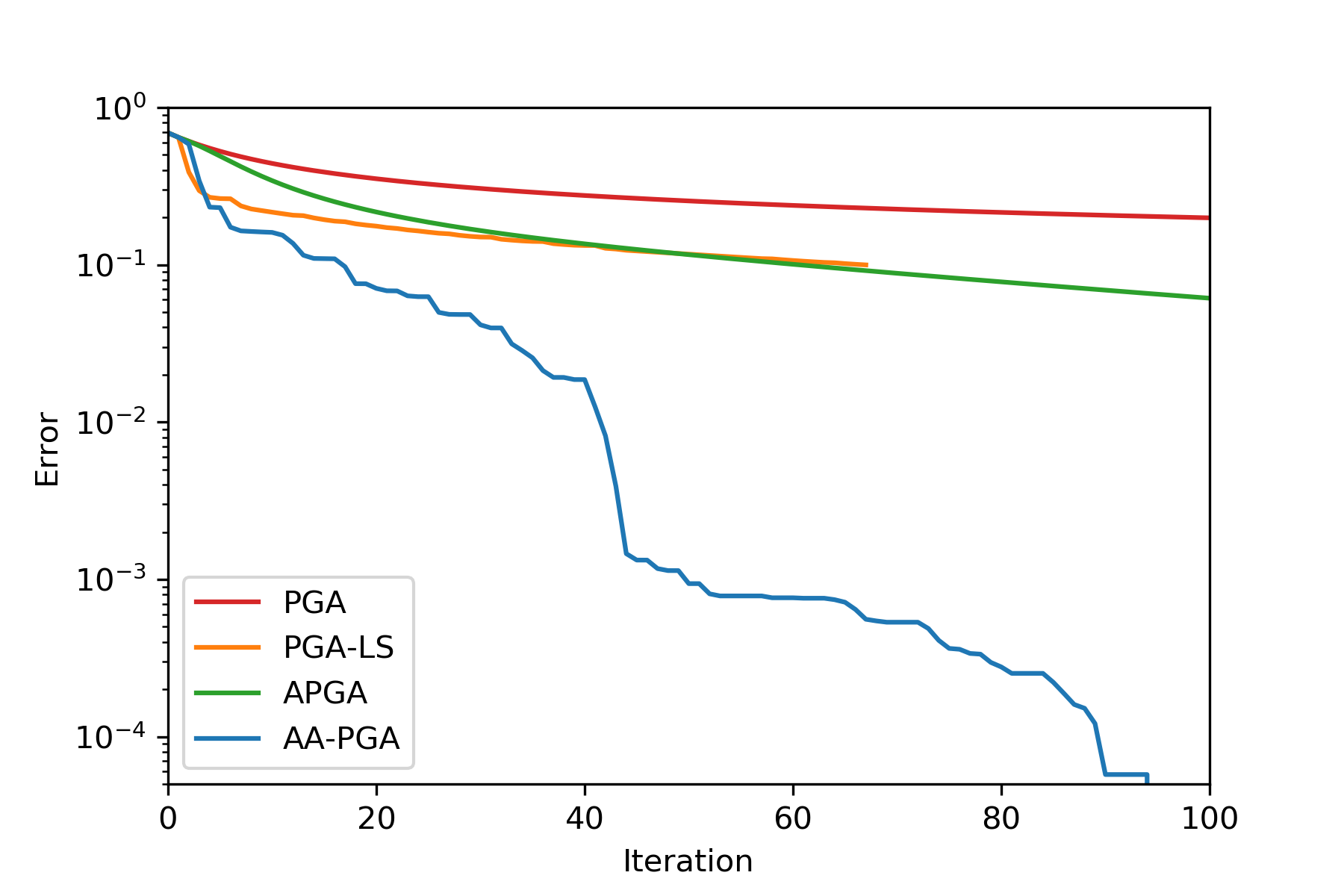}}
			{\includegraphics[width=0.4\textwidth]{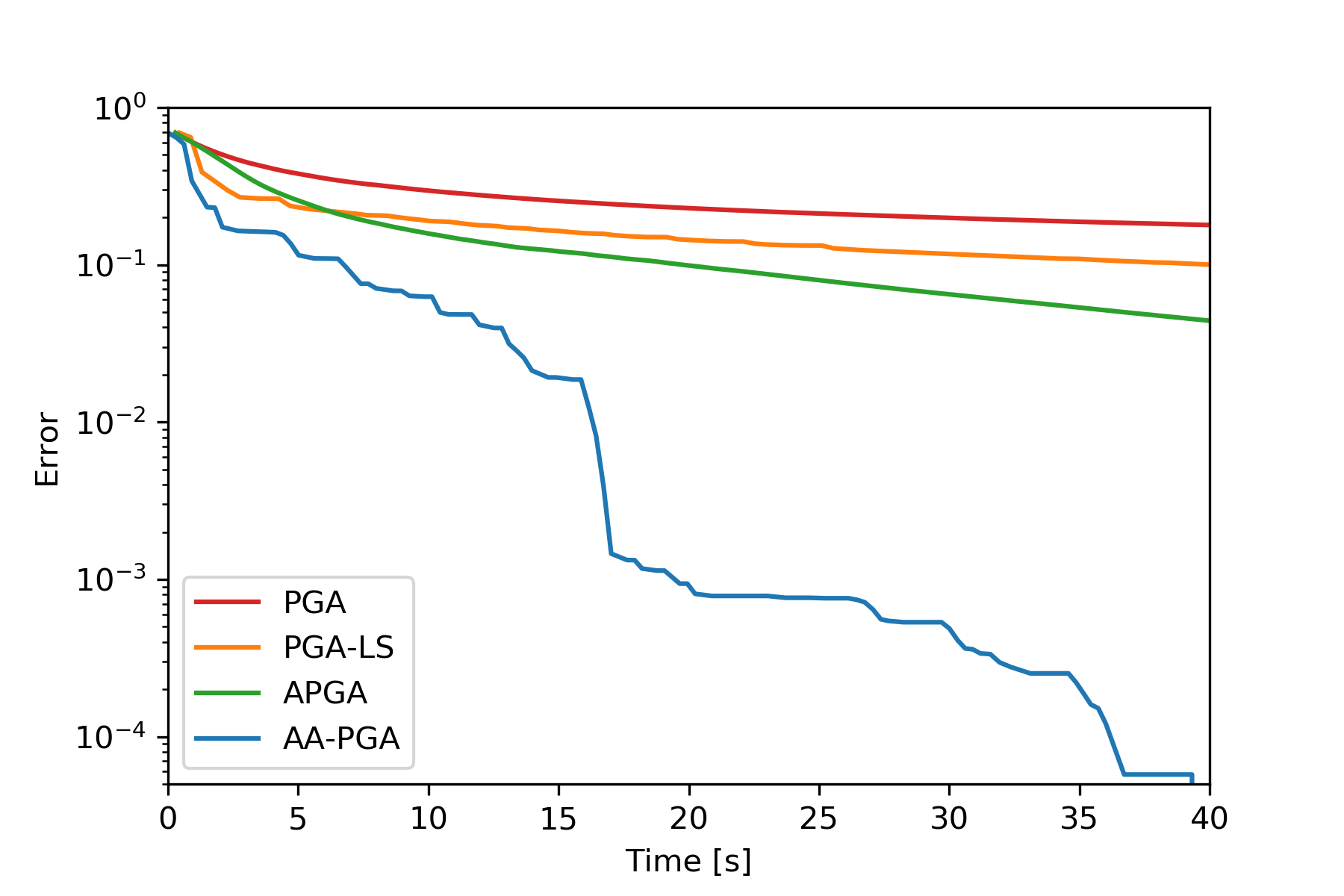}}
		}~\caption{Constrained logistic regression on the Madelon and the Gisette data sets.} \label{fig:logreg:1}
\end{figure}

\begin{figure}[t!]
	\centering 
		\subfigure[Cina0: $\mu=0.1$, $\kappa=1.2\times 10^{7}$]{
			{\includegraphics[width=0.4\textwidth]{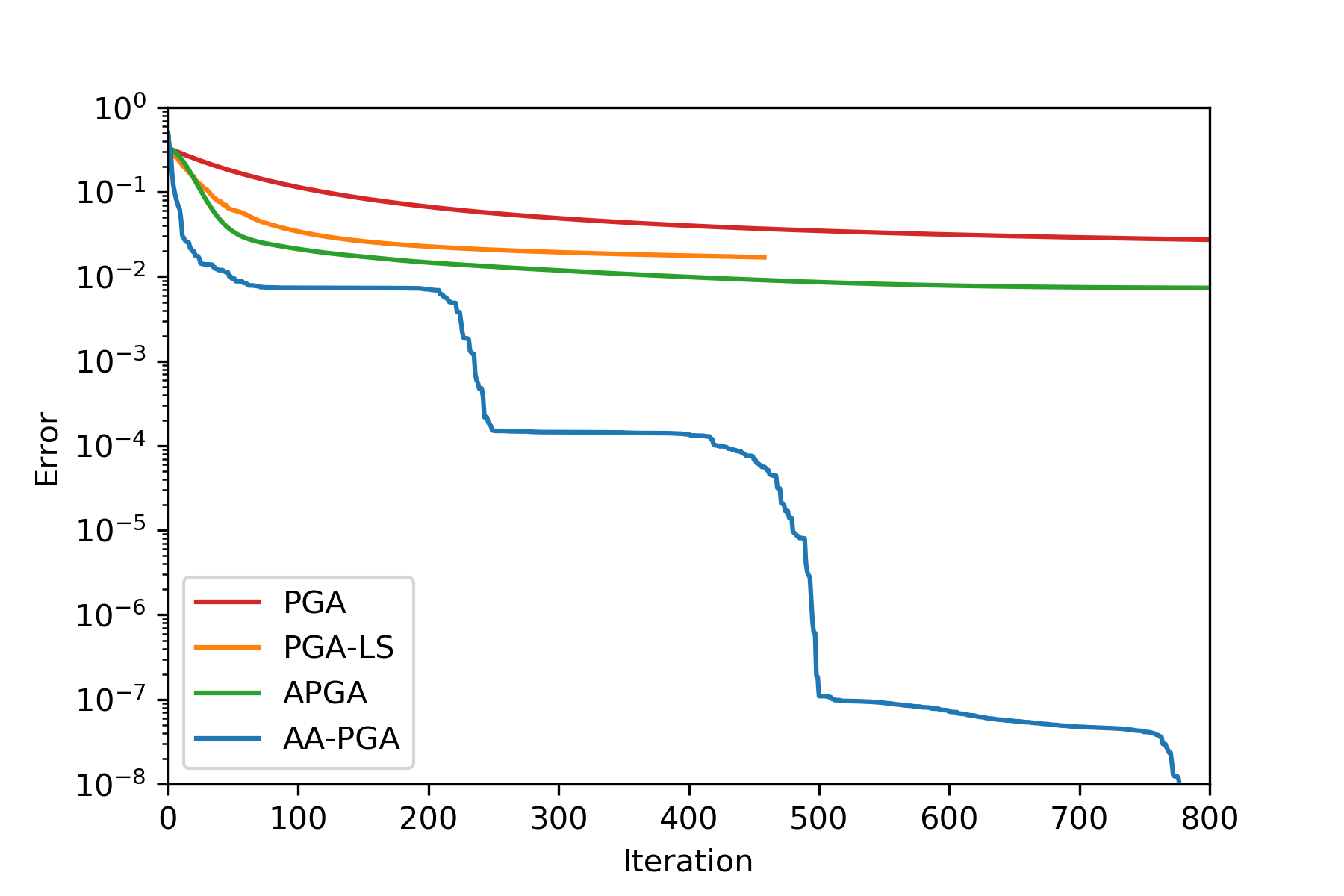}}
			{\includegraphics[width=0.4\textwidth]{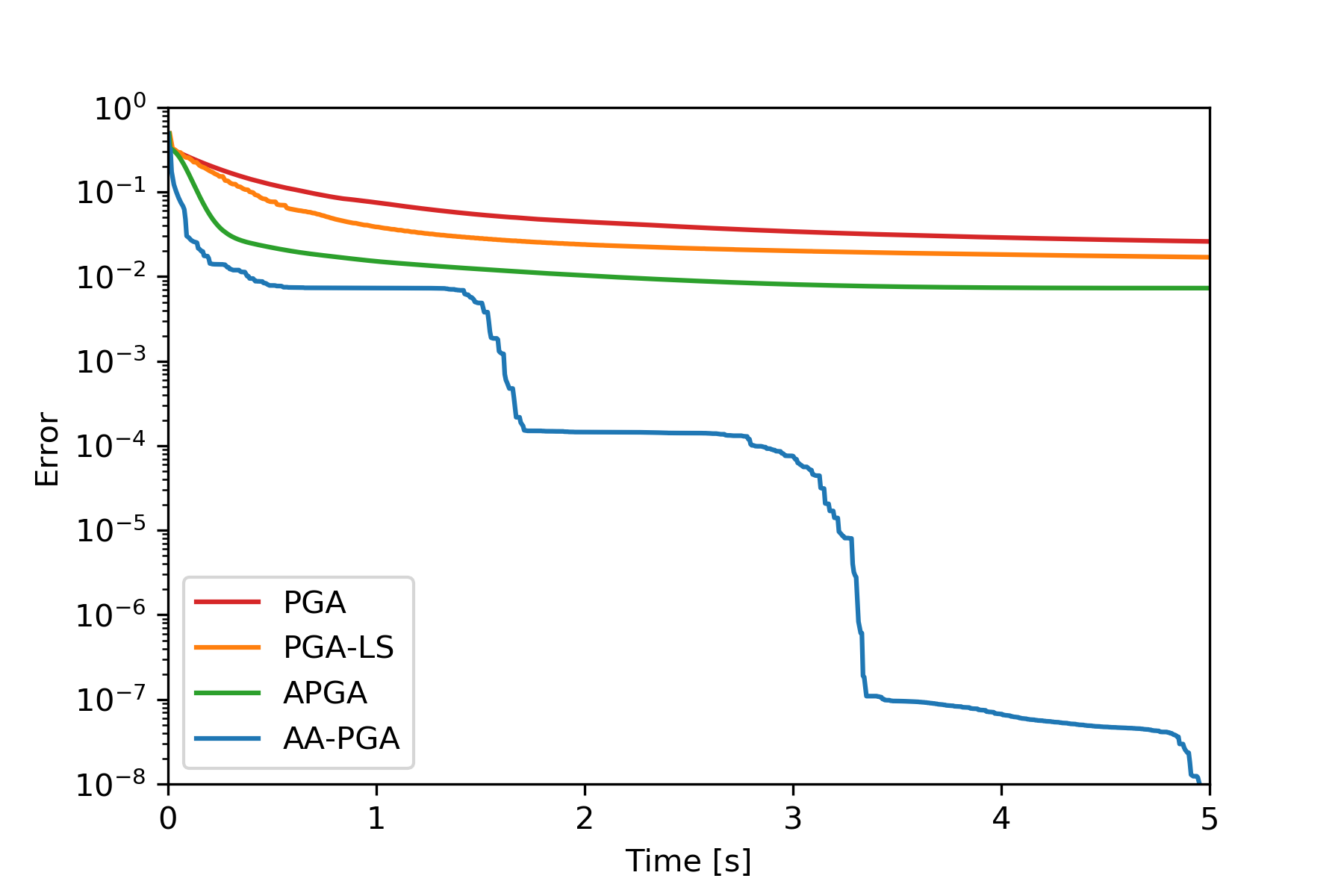}}
		}
 	\hfill
 	  \subfigure[Sido0: $\mu=10^{-2}$, $\kappa=3.7\times 10^{3}$]{
			{\includegraphics[width=0.4\textwidth]{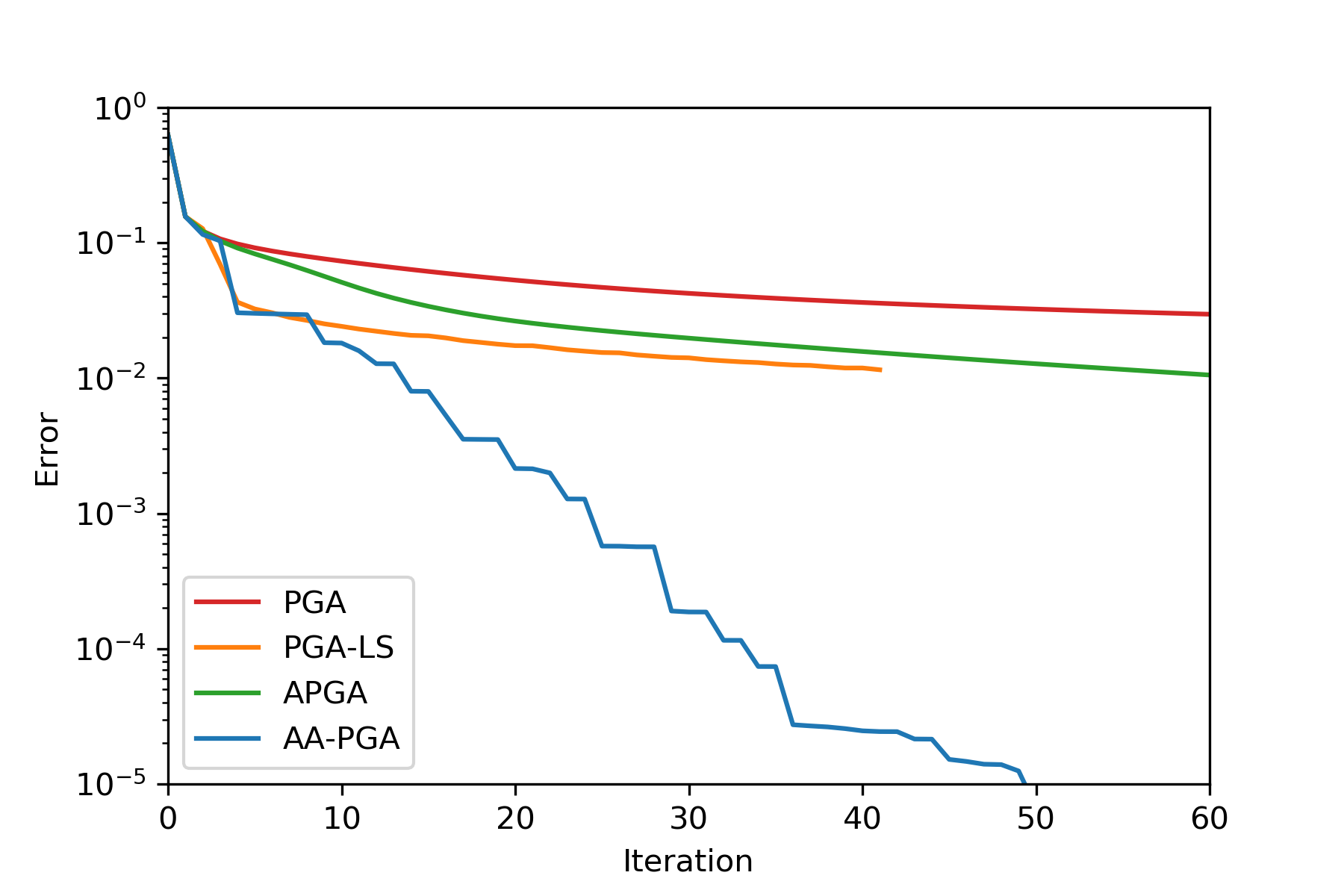}}
			{\includegraphics[width=0.4\textwidth]{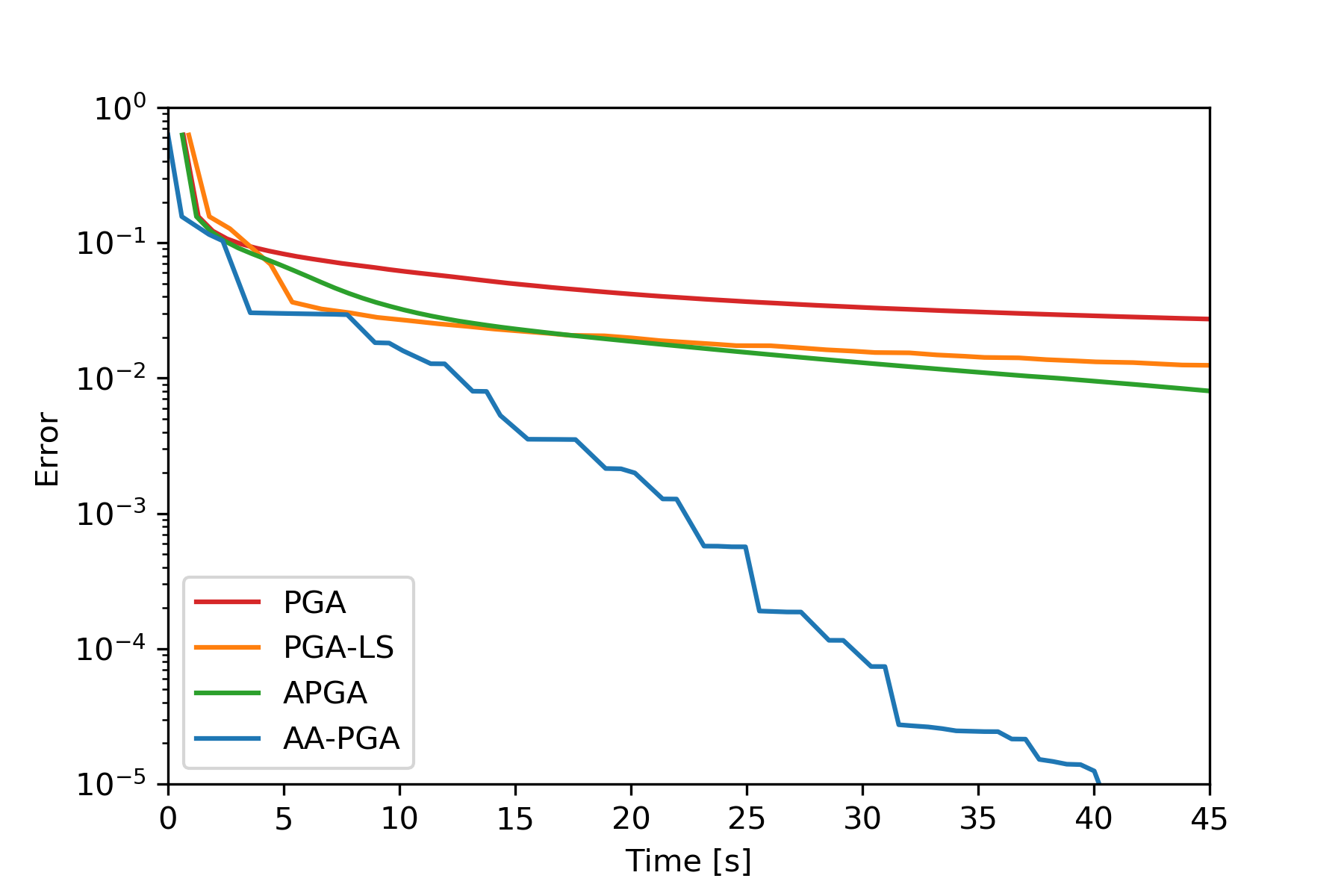}}
		}~\caption{Constrained logistic regression on the Cina0 and Sido0 data sets.} \label{fig:logreg:2}
\end{figure}

Figures~\ref{fig:logreg:1} and~\ref{fig:logreg:2} show the performance of AA-PGA and other selected algorithms on four different data sets. As can be seen, AA consistently and dramatically improves the performance of standard first order methods both in number of iterations  and wall-clock time. Since these data sets are very ill-conditioned, standard first order methods make very little progress, while AA can quickly find a high accuracy approximate solution. 
This once again demonstrates the great benefit of local adaptation and acceleration as previously seen in unconstrained quadratic problems (see., Figure~\ref{fig:quadratic}).  
In most cases, the convergence rate is linear confirming our prediction. The result also highlights the importance of the guard step in Algorithm~\ref{alg:aa:pga:gd}. Specifically, in some hard instances such as the one shown in  Fig.~\ref{fig:logreg:2}(a), the iterates alternate between periods with big jumps due to AA steps, which often significantly reduce the objective, and slowly converging regimes governed by the PGA steps. The later steps help to guide the iterates through a tough regime until AA steps take over and make big improvement.


\subsection{Nonnegative least squares} 

Next, we consider the nonnegative least squares problem: 
\begin{align*}
	\underset{x\in \R^n}{\mbox{minimize}} \,\frac{1}{2M}\norm{Ax-b}_2^2 + \mu \norm{x}_2^2 \quad \mbox{subject to}\,\,\,  x\geq 0,
\end{align*}
which is a core step in many nonnegative matrix factorization algorithms. We set $\gamma=1/L$, where $L=\norm{A}_2^2/M$. 

\begin{figure}[t!]
	\centering 
		\subfigure[Madelon: $\mu=0.1$, $\kappa=1.2\times 10^{9}$]{
			{\includegraphics[width=0.4\textwidth]{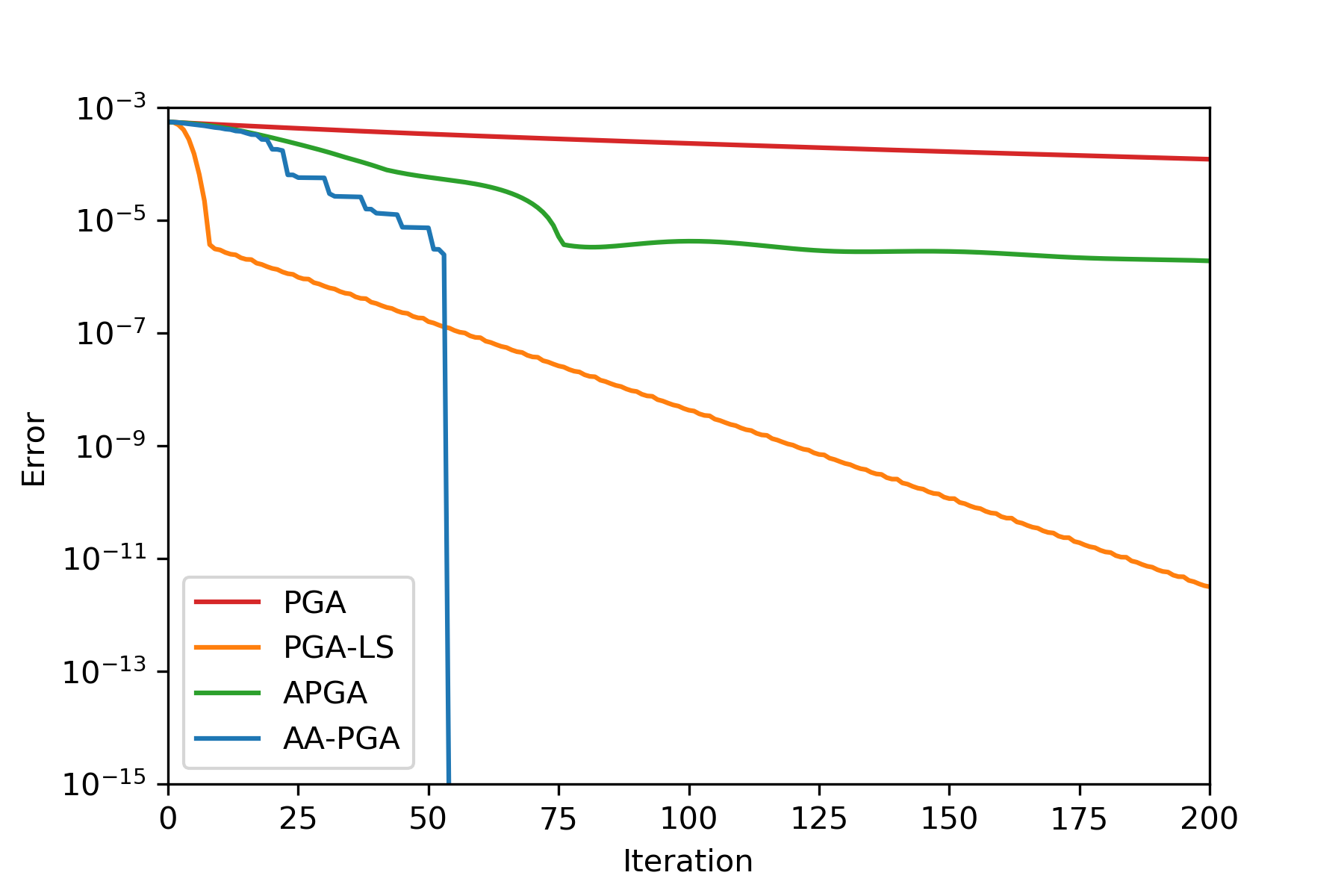}}
			{\includegraphics[width=0.4\textwidth]{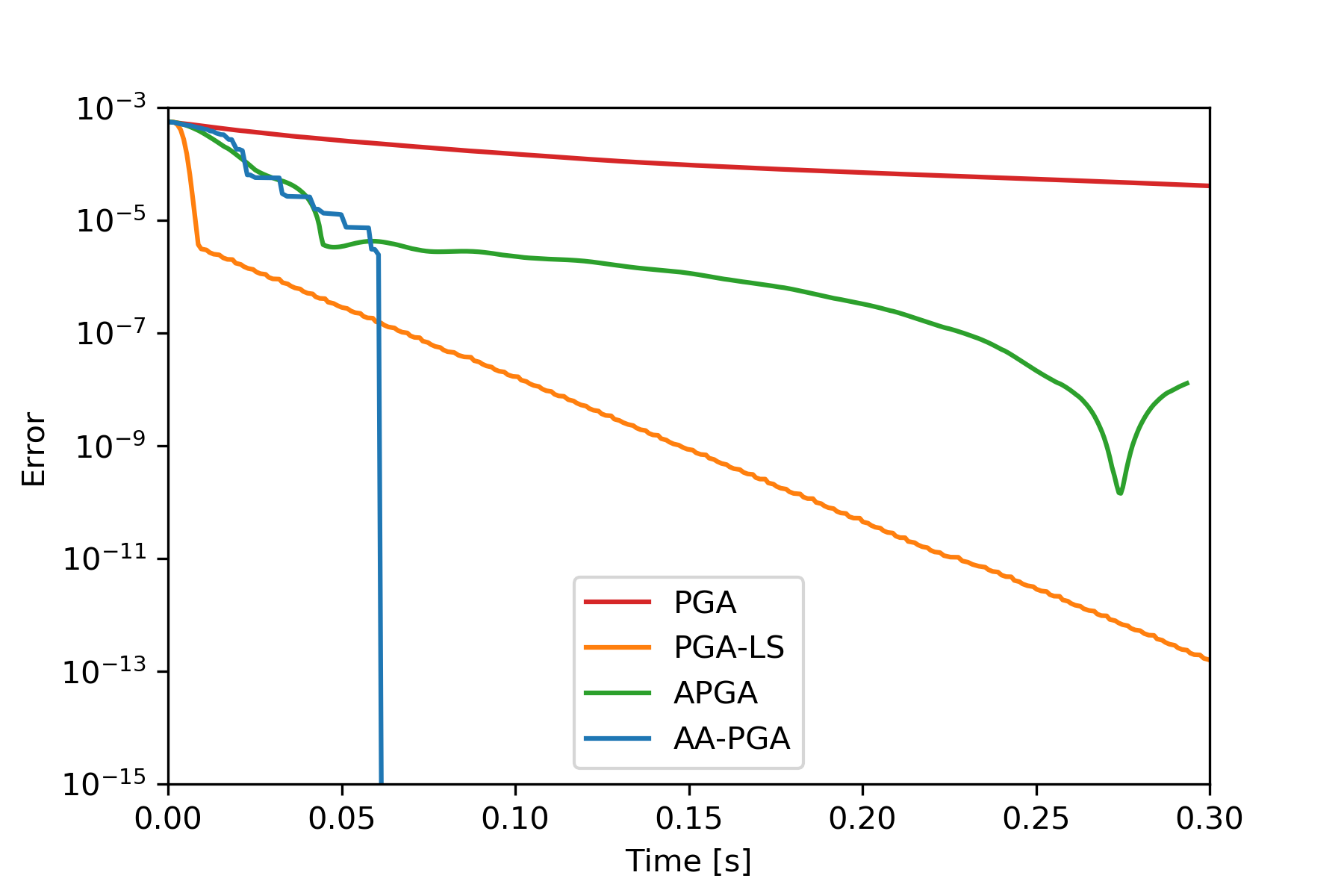}}
		}
	\hfill
		\subfigure[Gisette: $\mu=10$, $\kappa=1.36\times 10^{7}$]{
		{\includegraphics[width=0.4\textwidth]{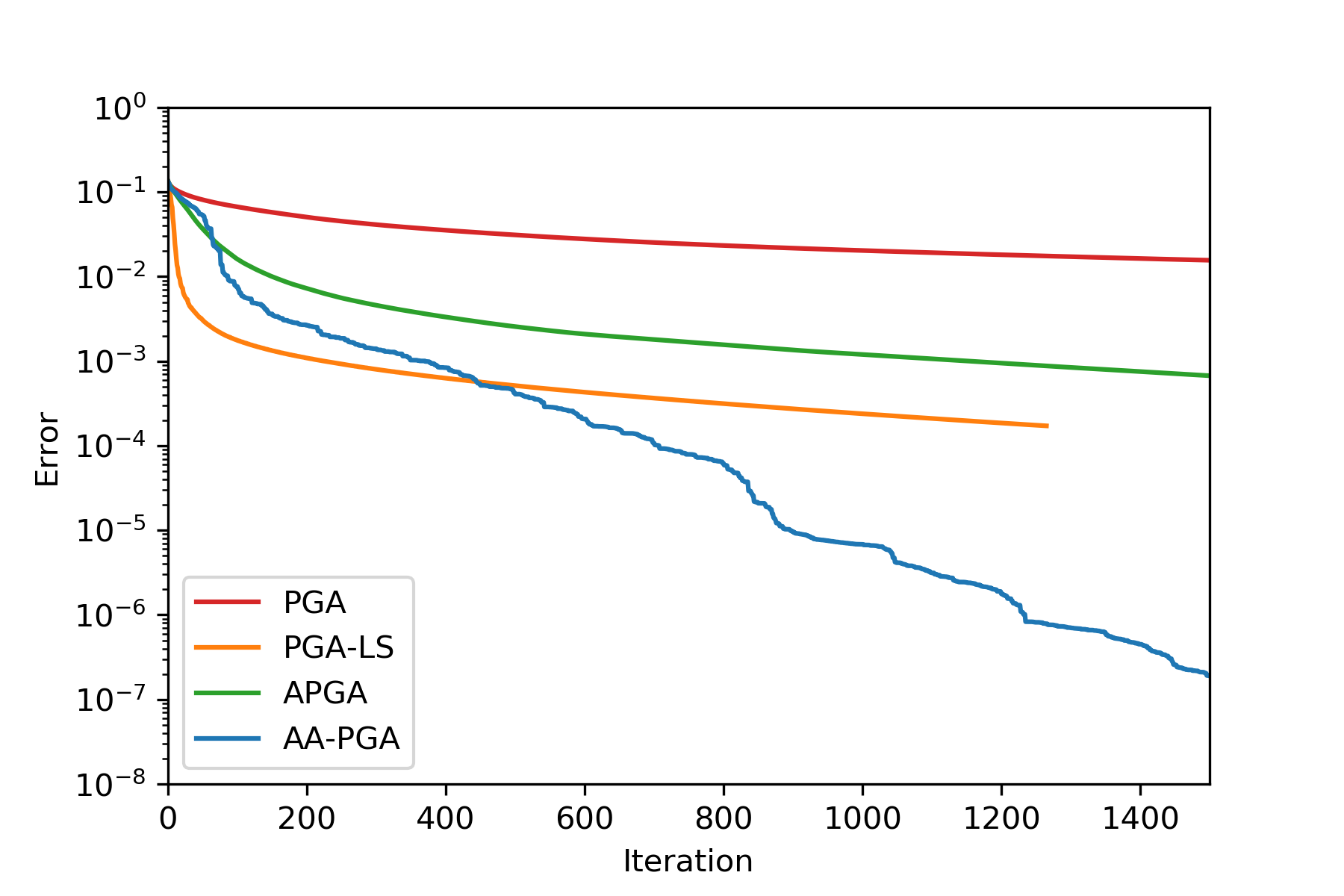}}
		{\includegraphics[width=0.4\textwidth]{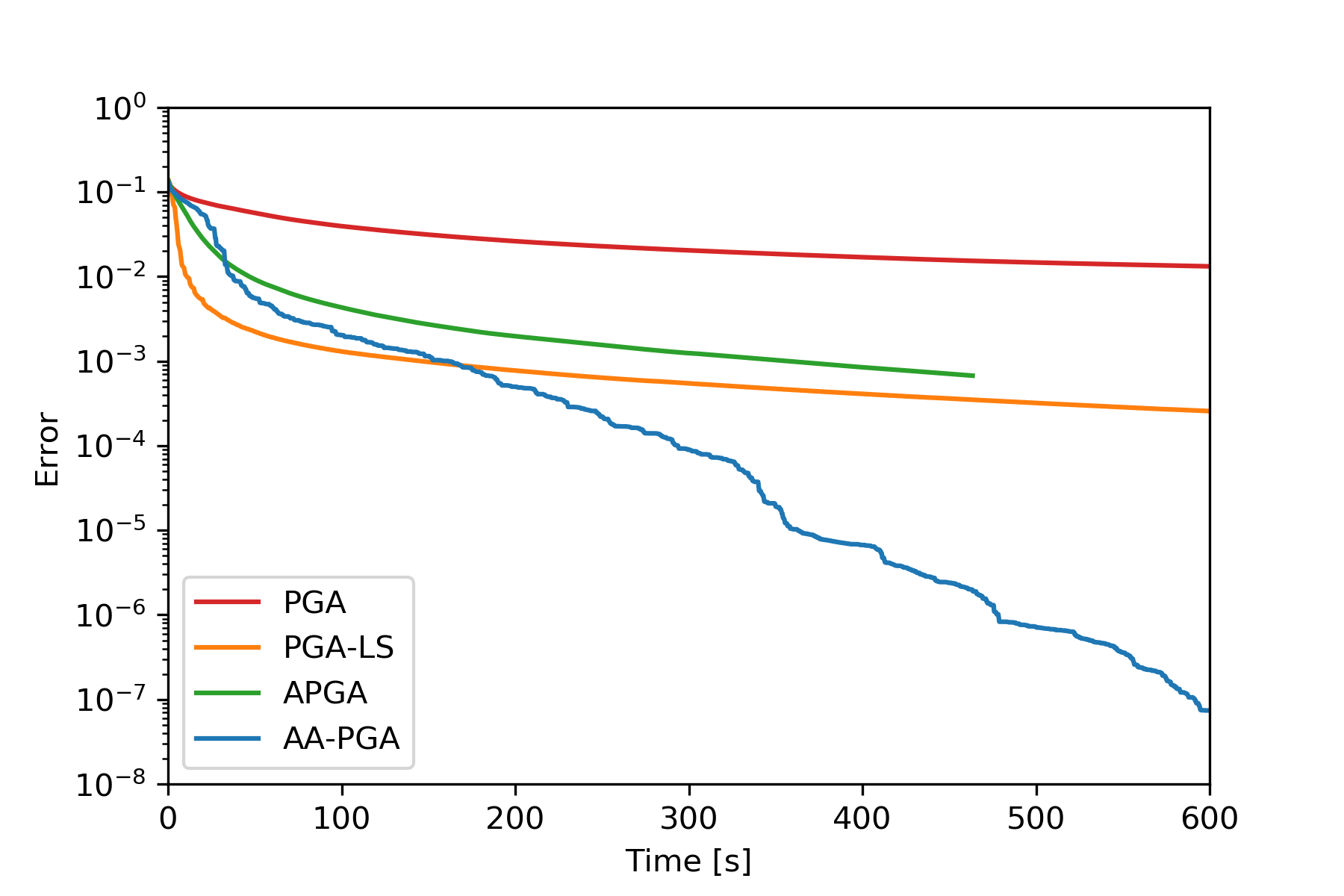}}		
	}~\caption{Nonnegative least-squares on the Madelon and Gisette data sets}\label{fig:nnls:1}
\end{figure}

Similarly to the previous problem, AA offers significant acceleration and often achieves several orders of magnitude speed-up over popular first order methods. Interestingly, in Fig.~\ref{fig:nnls:1}(a), AA seems to identify the solution in finite time. This could be the case where the optimal solution lies in the subspace spanned by the past iterates. 
\begin{figure}[t!]
	\centering 
		\subfigure[Cina0: $\mu=10$, $\kappa=4.8\times 10^{5}$]{
			{\includegraphics[width=0.4\textwidth]{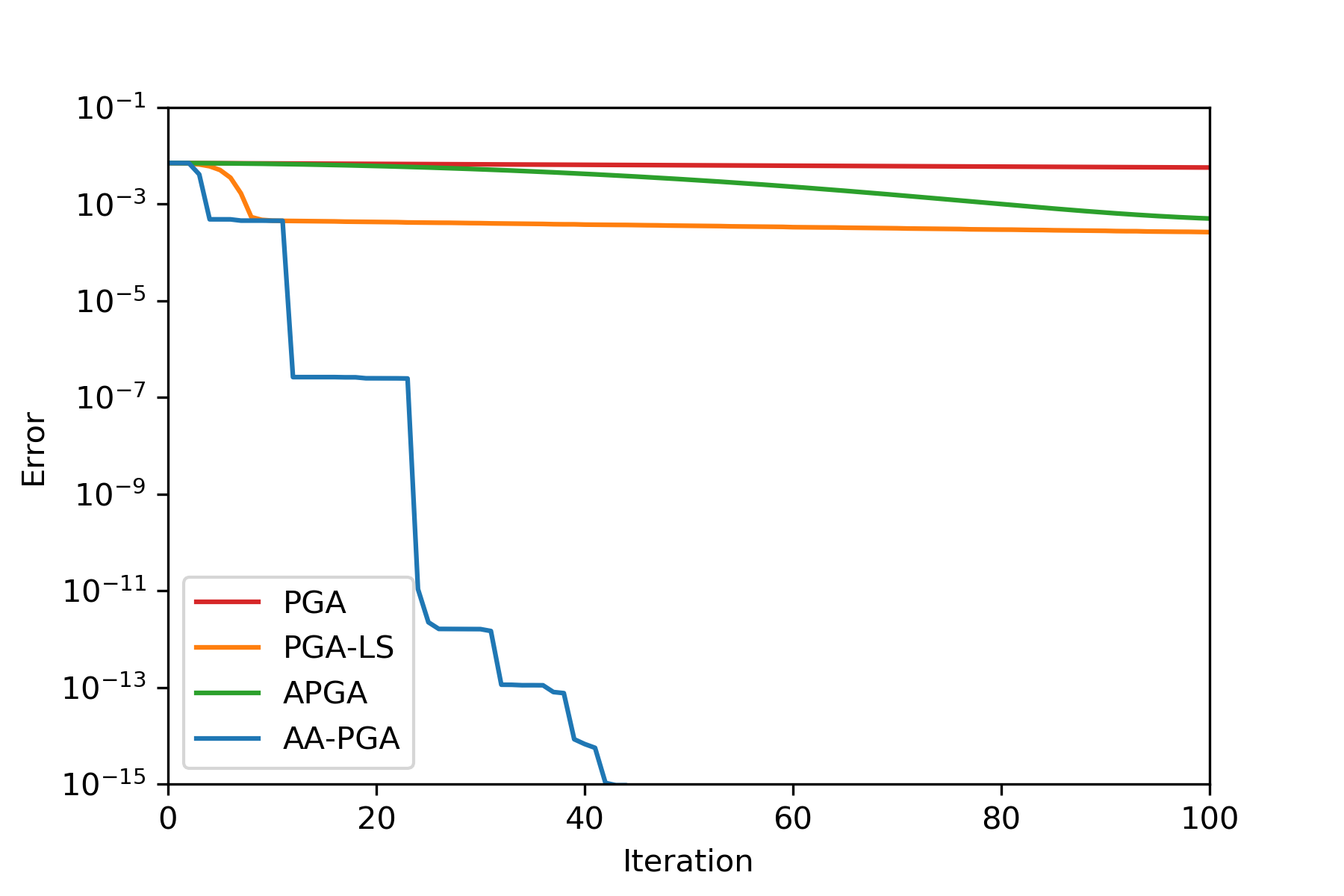}}
			{\includegraphics[width=0.4\textwidth]{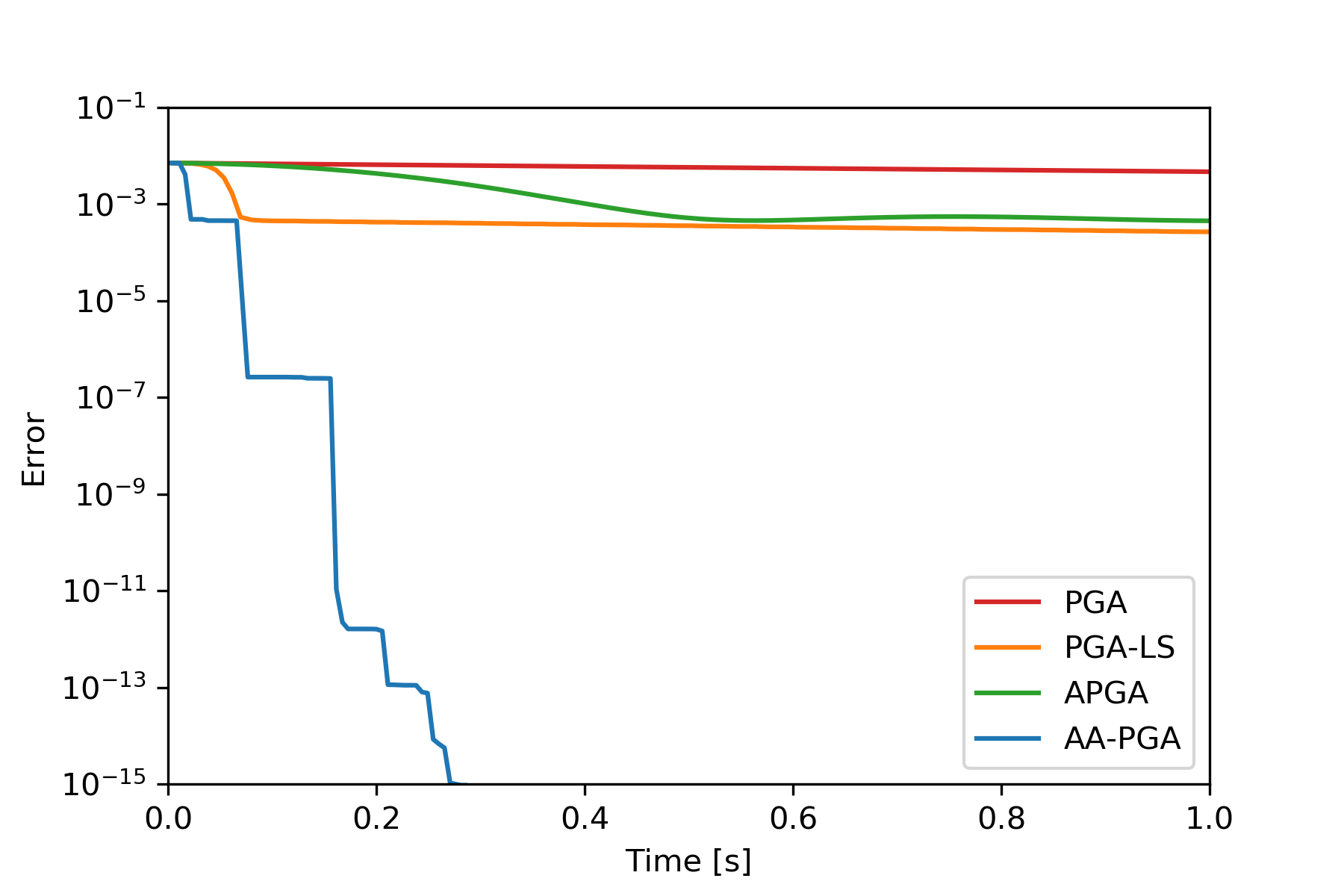}}			
		}
	\hfill
		\subfigure[Sido0: $\mu=0.1$, $\kappa=1.48\times 10^{6}$]{
			{\includegraphics[width=0.4\textwidth]{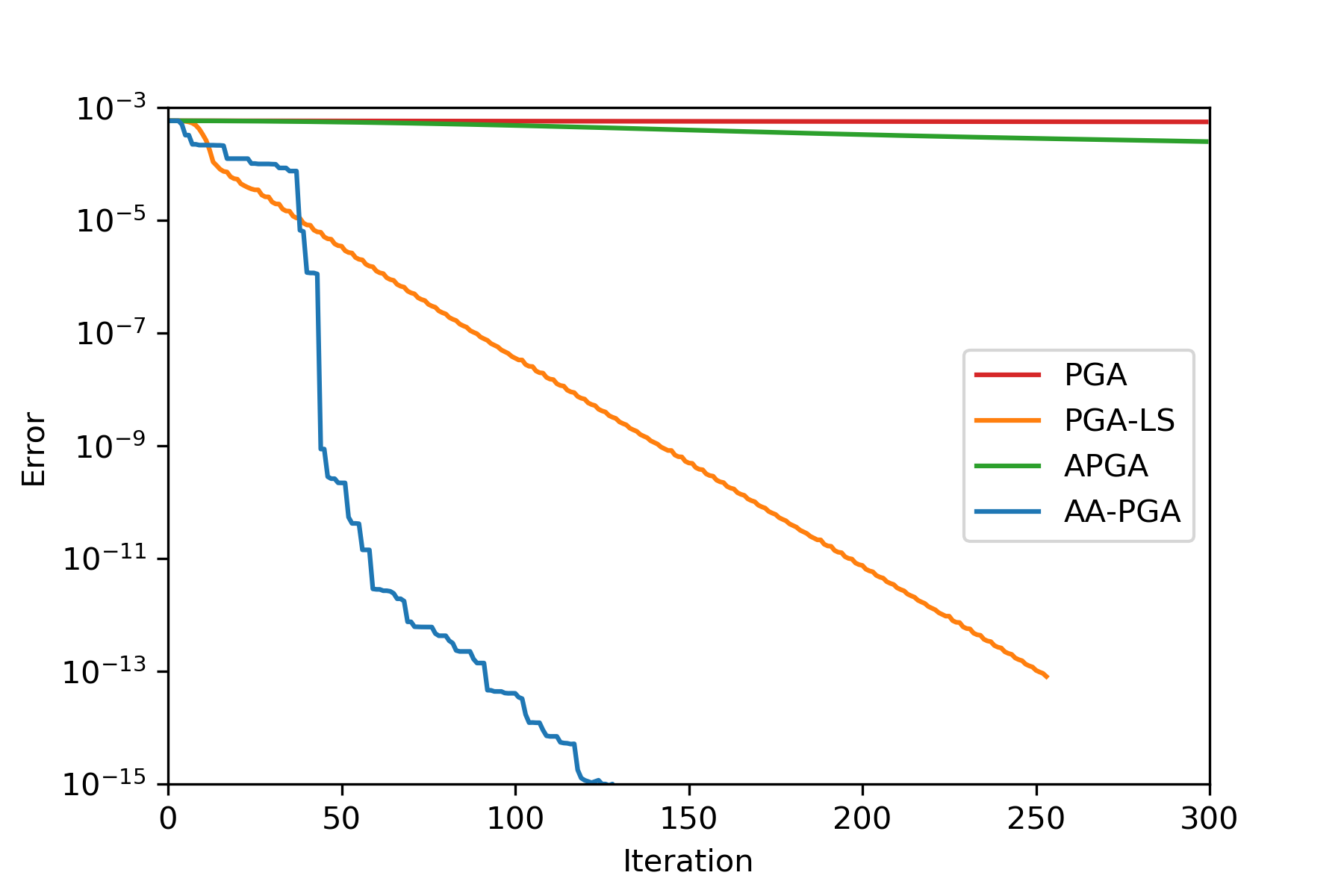}}
			{\includegraphics[width=0.4\textwidth]{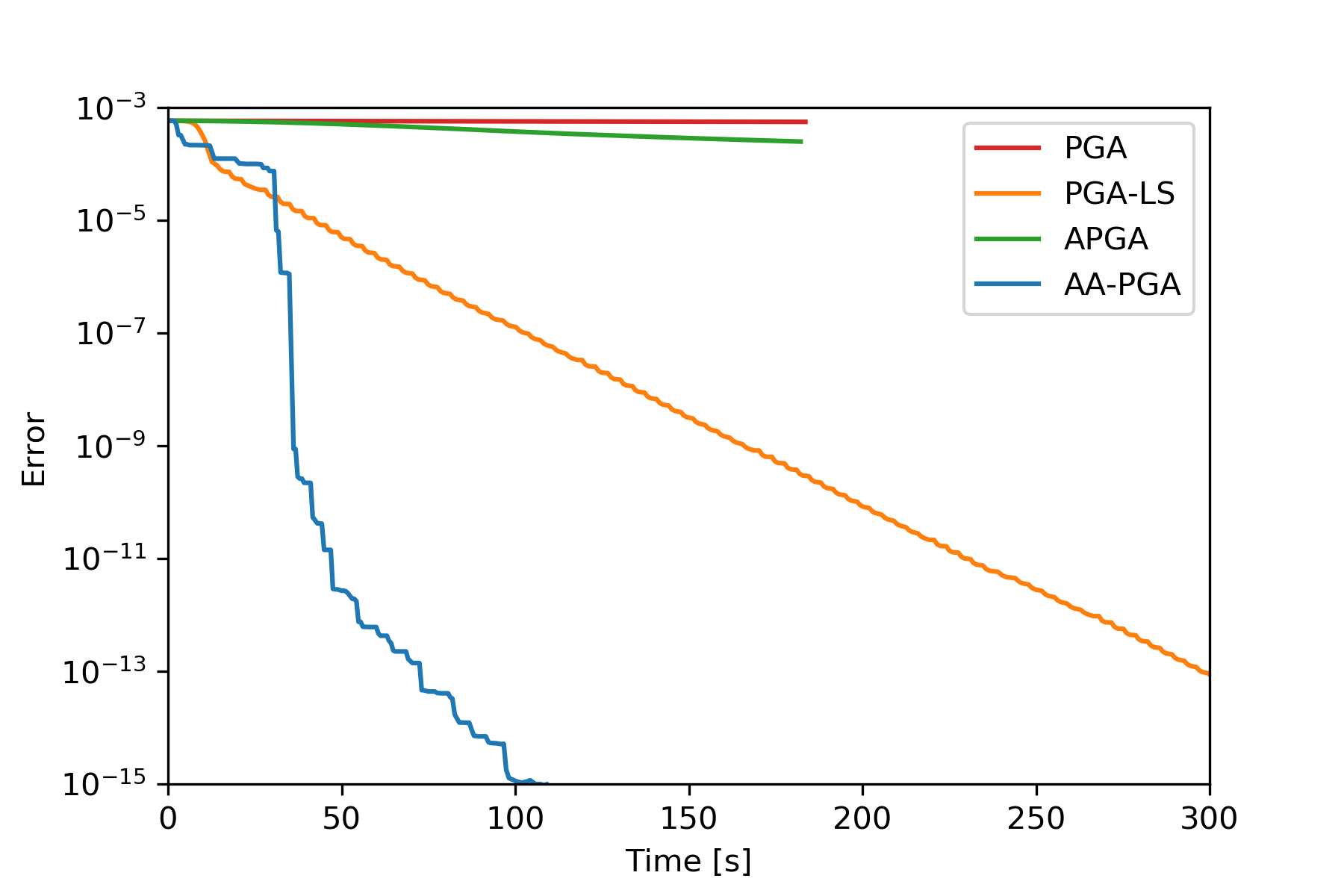}}		
		}~\caption{Nonnegative least-squares on the Cina0 and Sido0 data sets.} \label{fig:nnls:2}
\end{figure}

\begin{figure}[t!]
	\centering 
		\subfigure[$(m,n) = (100,1000)$]{
			{\includegraphics[width=0.4\textwidth]{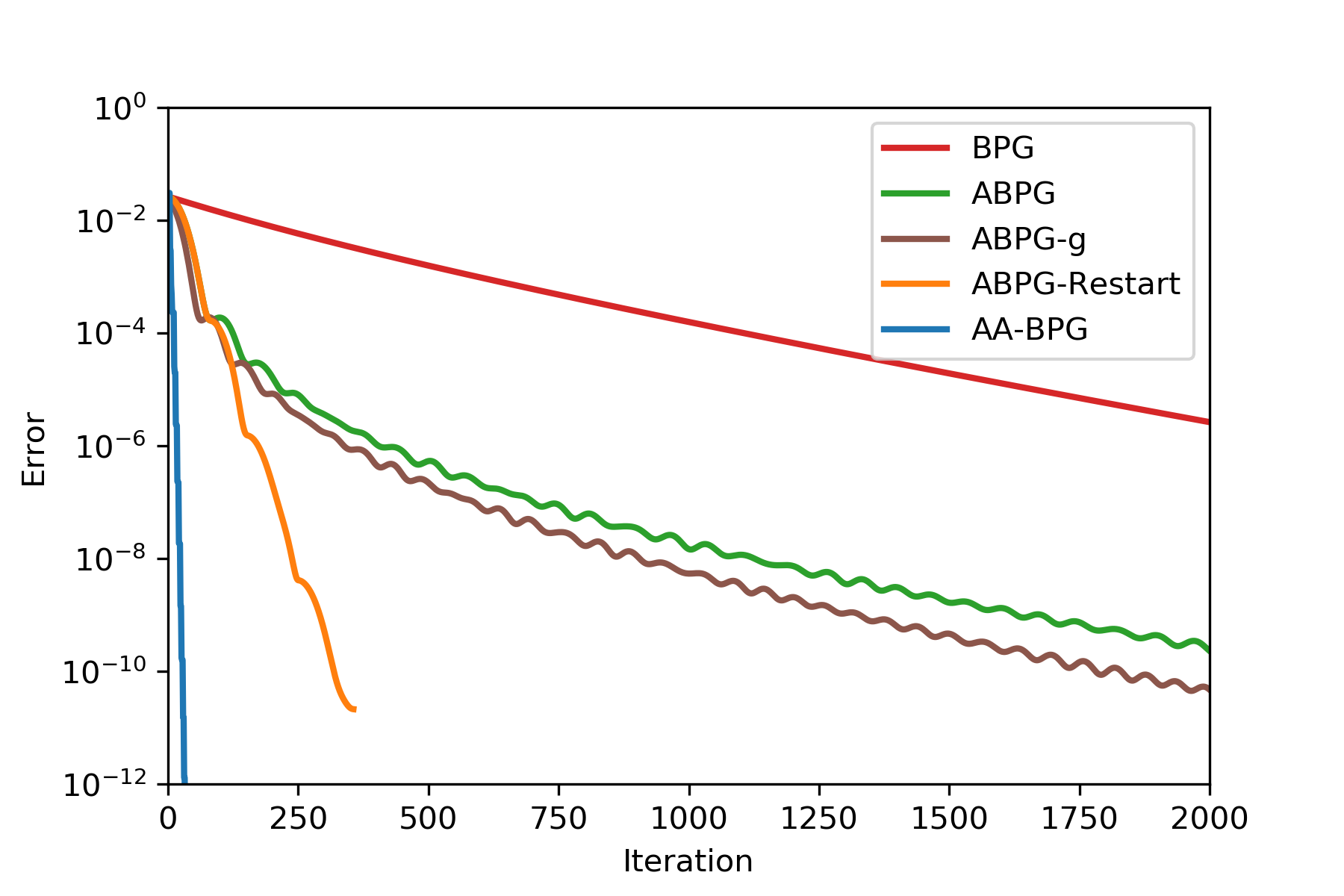}}
			\includegraphics[width=0.4\textwidth]{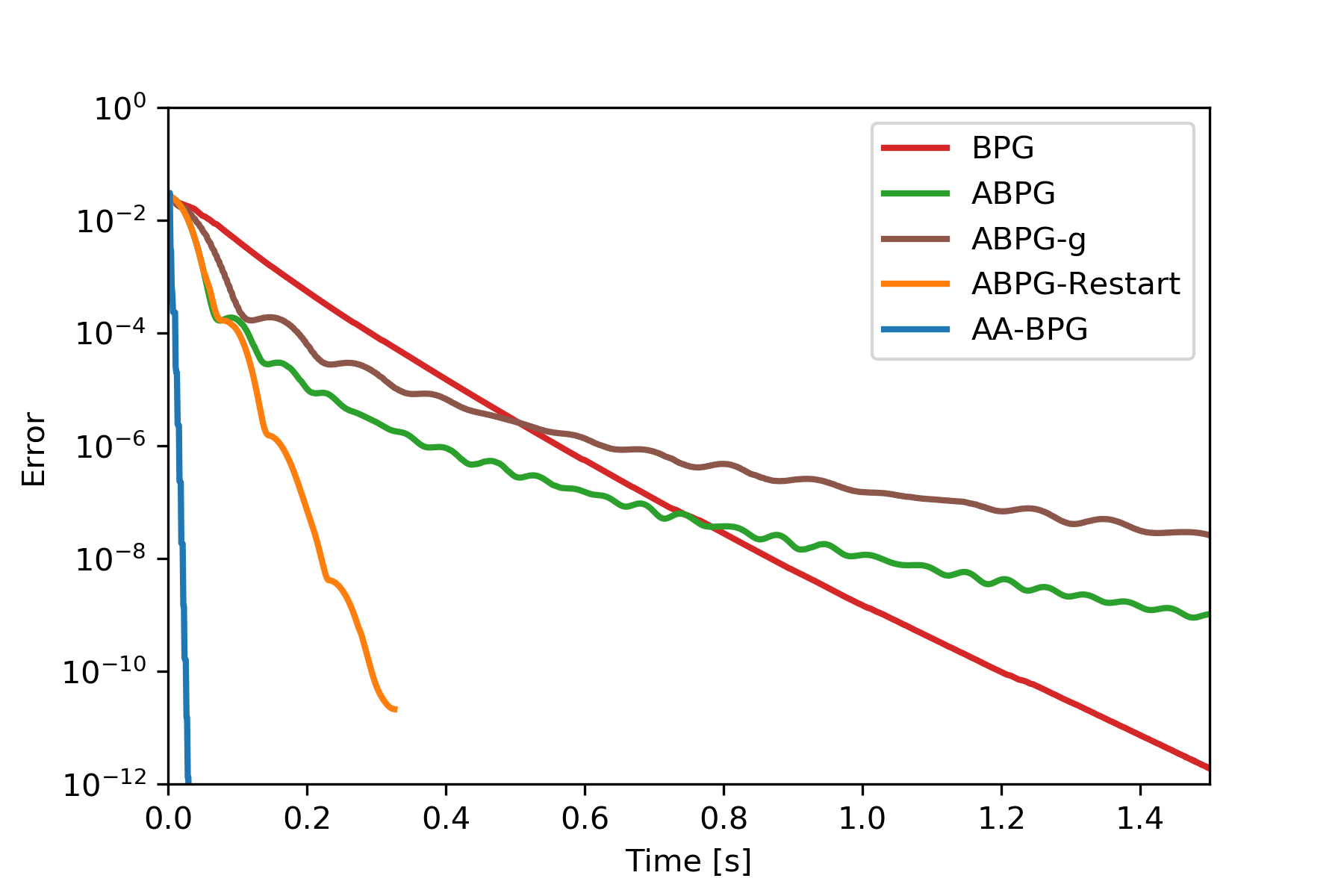}
		}
	\hfill 
	\subfigure[$(m,n) = (1000,100)$]{
		{\includegraphics[width=0.4\textwidth]{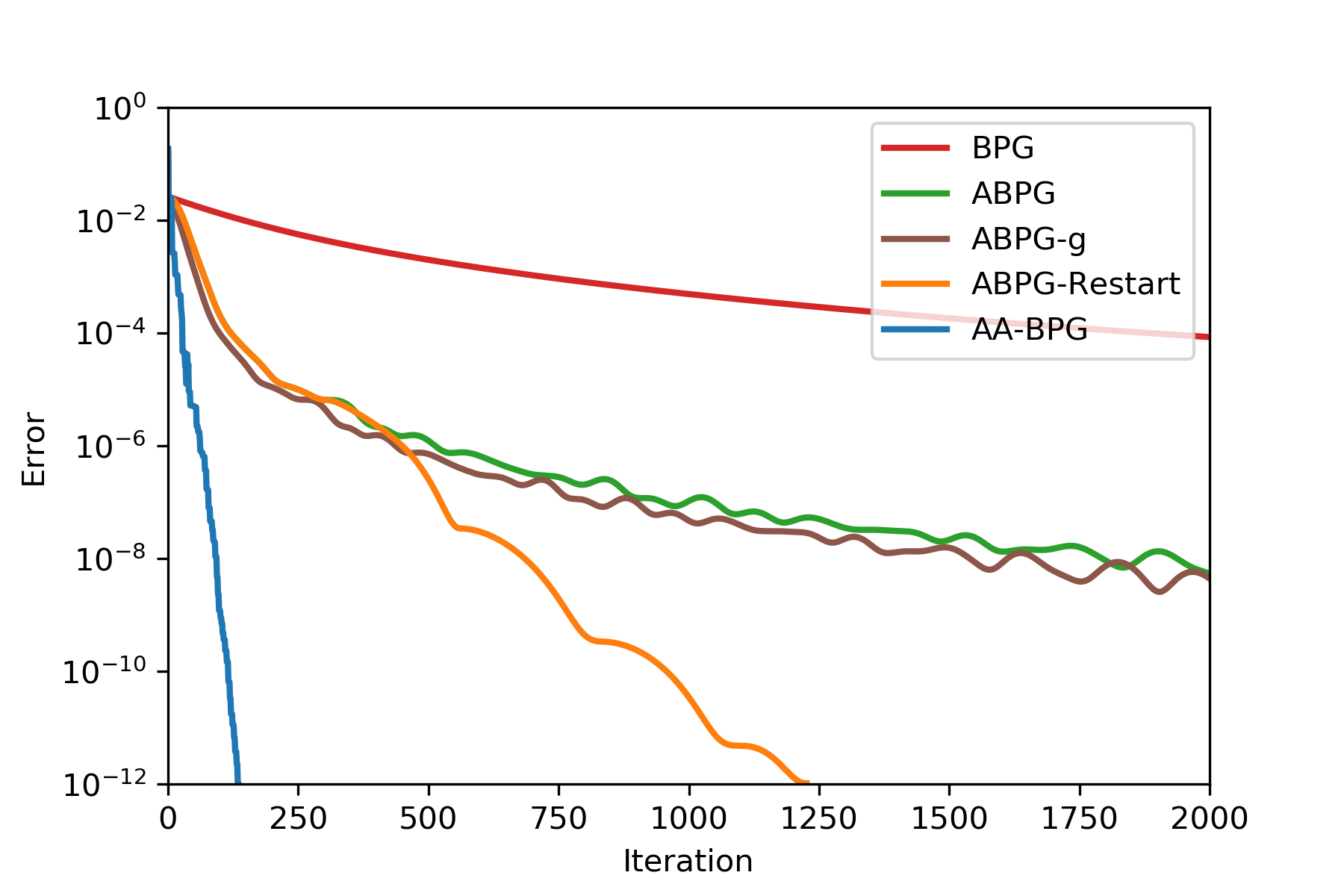}}
		\includegraphics[width=0.4\textwidth]{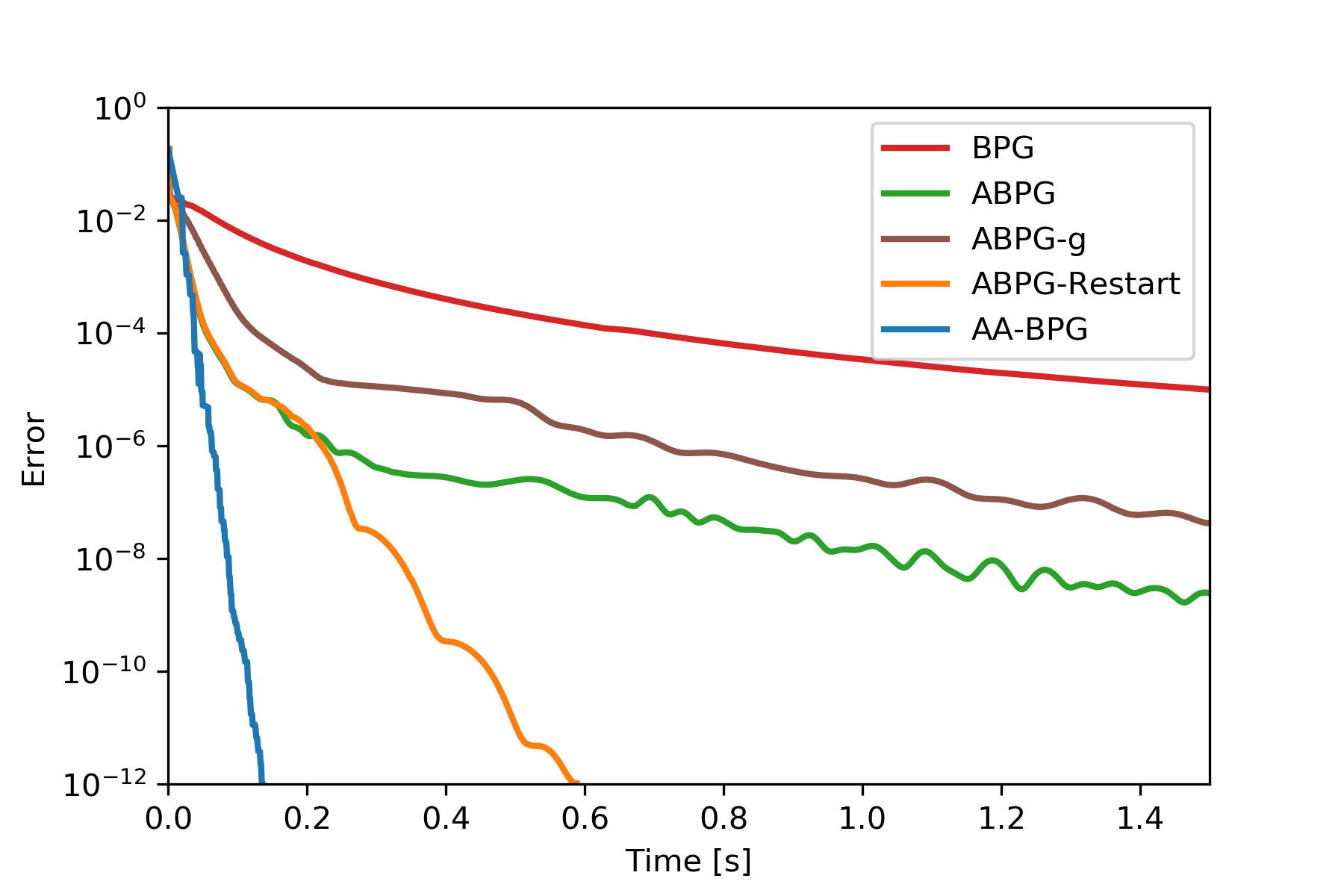}
	}~\caption{Relative-entropy nonnegative regression on two random problem instances.} \label{fig:poireg}
\end{figure}
\subsection{Relative-entropy nonnegative regression} The task is to reconstruct the signal $x\in\R^n_+$ by solving
\begin{align*}
	\underset{x}{\mbox{minimize}} \, D_{\mathrm{KL}}\left(Ax, b\right) + \lambda\norm{x}_1 \quad \mbox{subject to}\,\,\,  x\geq 0,
\end{align*}
where $A \in \R^{m \times n}_+$ is given nonnegative observation matrix and $b\in \R^m_{++}$ is a noisy measurement vector. 
We adapt the family of BPG methods with $\mathcal{D}=\R^n_+$, the Shannon entropy as the kernel $\genf$,  $f(x)=D_{\mathrm{KL}}\left(Ax, b\right)$, and $h(x) = \lambda \norm{x}_1$ with $\lambda=0.001$. It is shown in \cite{BBT16} that $f$ is $L$-smooth relative to $\genf$ with constant $L = \max_{1 \leq i \leq n} \norm{a_i}_1$. We follow \cite{HRX18} and generate two problem instances with $A$ and $b$ having entries uniformly distributed over the interval $[0,1]$. All methods are initialized at $x_0=\ones$.

Figure~\ref{fig:poireg}(a) shows the suboptimality for a randomly generated instance of the relative-entropy nonnegative regression problem with $m=100$ and $n=1000$. This instance is often referred as the easy case, and BPG converges linearly. Figure~\ref{fig:poireg}(b) shows similar results for the hard instance with $m=1000$ and $n=100$, where the BPG method converges sublinearly. In both cases, AA-BPG achieves the fastest convergence and significantly outperforms the others. Interestingly, AA-BPG is able to achieve linear convergence even in the hard case, which shows a clear evidence that our method adapts to the local strong convexity of the objective. This ability is observed consistently in all the problems and  data sets we have considered, and confirms our theoretical predictions. 

\section{Conclusion}

We adapted Anderson acceleration to proximal gradient methods, retaining their global (worst-case) convergence guarantees while adding the potential for local adaption and acceleration. Key innovations include theoretical convergence guarantees for non-smooth mappings, techniques for avoiding potential infeasibilities, and stabilized algorithms with global convergence rate guarantees and strong practical performance.
We also proposed an application of AA to non-Euclidean geometry. Given that AA can be applied to general fixed-point computations, the current literature has just scratched the surface of potential uses of AA in optimization. With its simplicity and evident promise, we feel that AA merits much further study.

\section*{Acknowledgements}

This work was supported in part by the Knut and Alice Wallenberg Foundation, the Swedish Research Council and the Swedish Foundation for Strategic Research. We would like to thank Wenqing Ouyang for his useful feedback and suggestions on an early draft of this paper. We also thank the anonymous reviewers for their useful comments and suggestions.

\bibliographystyle{abbrv}
\bibliography{refs}

\appendix

\section{Proof of Theorem~\ref{thrm:main}}\label{appendix:thrm:main}

Since $g$ is differentiable at $y\opt$ with Jacobian $G$,  it holds  that \cite[Eq. 9(6)]{RW07}:
\begin{align}\label{eq:gy:linearization}
	g(y) = y^\star + G(y-y^\star) + e(y),
\end{align}
where $e(y) = o\left(\ltwo{y-y^\star}\right)$ as $\ltwo{y-y^\star} \to 0$. This means that for any $\varepsilon>0$, there exits $\delta>0$ such that
\begin{align}\label{eq:gy:linearization:norm}
%	\ltwo{e(y)} \leq \varepsilon \ltwo{y-y\opt}, \quad \forall y \in \mc{B}(y\opt, \delta),
	\ltwo{g(y) - y^\star - G(y-y^\star)} \leq \varepsilon \ltwo{y-y\opt} &&\forall y \in \mc{B}(y\opt, \delta),
%\end{align} 
%We thus have for $y \in \mc{B}(y\opt, \delta)$ that
\intertext{and hence}
%\begin{align*}
	\ltwo{g(y) - y^\star} \leq \left({\ltwo{G}+\varepsilon}\right) \ltwo{y-y\opt} &&\forall y \in \mc{B}(y\opt, \delta)\nonumber.
\end{align}
Take $\varepsilon < (1- \ltwo{G})/2$ and define $\rho = \ltwo{G} + \varepsilon$. 
Note that   $\rho \in(0,1)$ and that $\rho+\varepsilon<1$.
Since $F(y) = g(y)-y$ with $F(y^\star)=0$, it holds for any $y\in\mc{B}(y\opt, \delta)$ that:
\begin{align*}
\norm{F(y)-F(y\opt)}_2
	\leq
		\norm{g(y)-y\opt}_2
		+
		\norm{y-y\opt}_2
%	\\
	\leq
	\left(1+\rho\right)\norm{y-y\opt}_2.
\end{align*}
Similarly, we have 
\begin{align*}
\norm{y-y\opt}_2
	\leq
		\norm{g(y)-y\opt}_2
		+
		\norm{F(y)-F(y\opt)}_2
%	\\
	\leq
	\rho \norm{y-y\opt}_2
	+
	\norm{F(y) -F(y\opt)}_2.
\end{align*}
In summary, it holds for any $y\in\mc{B}(y\opt, \delta)$ that
\begin{align}\label{eq:Fy:twoside}
	 {(1-\rho)}\norm{y-y^\star}_2 
	 \leq 
	 \norm{F(y)}_2 
	 \leq {(1+\rho)}\norm{y-y^\star}_2.
\end{align}

Now, let $\hat{\rho} = \rho + \varepsilon'$ for some $\varepsilon'>0$ such that $\hat{\rho}\in (\rho,1)$. We will show by induction that when $y_0$ is sufficiently close to $y\opt$, we have
\begin{align*}
	\ltwo{F(y_k)} \leq \hat{\rho}^k \ltwo{F(y_0)} \quad \forall k.	
\end{align*}
To that end, we will pick an $\varepsilon < (1 -\ltwo{G})/2$ above small enough such that
\begin{align}\label{eq:y:init:cond:r}
		\frac{\rho}{\hat{\rho}}
		+
		\frac{\varepsilon M_\alpha \,\hat{\rho}^{-m-1}}{1-\rho}
		\leq 
			1- \frac{\varepsilon}{1-\rho}.
\end{align}
Let $\delta$ be determined by the chosen $\varepsilon$. Fix a radius $r\leq \delta$ and let $y_0\in\mc{B}(y\opt, r)$ satisfy
\begin{align}\label{eq:y:init:cond:r0}
		\frac{M_\alpha (1+\rho) }{1-\rho} \hat{\rho}^{-m} \norm{y_0-y\opt}_2 
		\leq 
			r \qquad \forall y_0\in\mc{B}(y\opt,r).
\end{align}
We can now proceed as~\cite{TK15}. First,  note that the base case $k=0$ is obvious. Next, suppose that the hypothesis is true up to iteration $k$. We deduce for all $i=0,1, \ldots, m_k$ that:  
\begin{align}\label{eq:hypothesis:y}
	\ltwo{y_{k-i} - y\opt} \leq \frac{\hat{\rho}^{k-i}}{1- \rho}  \ltwo{F(y_0)} 
	\leq
		\frac{1+\rho}{1- \rho} \hat{\rho}^{k-i} \ltwo{y_0-y\opt},
\end{align}
where we used the induction hypothesis and \eqref{eq:Fy:twoside}.
It follows from \eqref{eq:y:init:cond:r0} that $y_{k-i}\in\mc{B}(y\opt, \delta)$ for all $i=0,1, \ldots, m_k$, and hence by \eqref{eq:gy:linearization}, we have
\begin{align}\label{eq:gy:linearization:k}
	g(y_{k-i}) = y\opt + G(y_{k-i}-y\opt) + e_{k-i},
\end{align}
where $e_{k-i}:=e(y_{k-i})\leq \varepsilon \ltwo{y_{k-1}-y\opt}$.
Thus, $y_{k+1} = \sum_{i=0}^{m_k} \alpha_i^k g(y_{k-i})$ can be written as 
\begin{align}\label{eq:ynext}
	y_{k+1} = y\opt + \sum_{i=0}^{m_k} \alpha_i^k\left[ G(y_{k-i}-y\opt) +  e_{k-i}\right].
\end{align}
Let $\bar{e}_k = \sum_{i=0}^{m_k} \alpha_i^k e_{k-i}$, we have 
\begin{align}\label{eq:ynext:term:1}
	\ltwo{\bar{e}_k} 
	\leq 
		\sum_{i=0}^{m_k} |\alpha_i^k| \ltwo{e_{k-i}}
	&\mathop \leq \limits^{\mathrm{(a)}} 
		\sum_{i=0}^{m_k} |\alpha_i^k| \varepsilon \ltwo{y_{k-i}-y\opt}
	\nonumber\\
	&\mathop \leq \limits^{\mathrm{(b)}} 
	 	\frac{\varepsilon M_\alpha}{1-\rho}\hat{\rho}^{k-m}\ltwo{F(y_0)}
	\nonumber\\
	&\leq
	 	\frac{\varepsilon M_\alpha}{1-\rho}\hat{\rho}^{-m}\ltwo{F(y_0)}, 
\end{align}
where $(\mathrm{a})$ follows from \eqref{eq:gy:linearization:norm}, and $(\mathrm{b})$ follows from Assumption~\ref{assumption:bounded:extra:coeffs}, the first inequality in \eqref{eq:hypothesis:y} and the fact that $m_k\leq m$.
We also have
\begin{align}\label{eq:ynext:term:2}
	\ltwo{\sum_{i=0}^{m_k} \alpha_i^k  G(y_{k-i}-y\opt) }
	\leq
		\sum_{i=0}^{m_k} |\alpha_i^k| \ltwo{G} \ltwo{y_{k-i}-y\opt}
	\leq 
		 \frac{\rho M_\alpha }{1-\rho} \hat{\rho}^{-m}\ltwo{F(y_0)},
\end{align}
where we used \eqref{eq:hypothesis:y} and $\ltwo{G}\leq \rho$ in the last step.
Combining \eqref{eq:ynext}, \eqref{eq:ynext:term:1}, and \eqref{eq:ynext:term:2} yields
\begin{align*}
	\ltwo{y_{k+1} - y\opt} 
	\leq
		\frac{M_\alpha \left(\rho  + \varepsilon \right)}{1-\rho} 
		\hat{\rho}^{-m}\ltwo{F(y_0)}.
\end{align*}
In view of \eqref{eq:Fy:twoside} and \eqref{eq:y:init:cond:r0}, it holds that $\ltwo{y_{k+1} - y\opt} \leq r$.  Hence,
%We can then express $g(y_{k+1})$ as
\begin{align*}
	g(y_{k+1}) = y\opt + G(y_{k+1} - y\opt) + e_{k+1},
\end{align*}
where $e_{k+1}$ satisfies
\begin{align}\label{eq:enext}
	\ltwo{e_{k+1}} \leq \varepsilon \ltwo{y_{k+1}-y\opt} \leq \varepsilon/(1-\rho)\ltwo{F(y_{k+1})}.
\end{align}
Since $F(y_{k+1}) = (G-\IM{})(y_{k+1} - y\opt) + e_{k+1}$,
it follows from \eqref{eq:ynext}  that
\begin{align}\label{eq:Fnext:expansion}
 	F(y_{k+1}) 
 	&= 
 		(G-\IM{})\sum_{i=0}^{m_k} \alpha_i^k G(y_{k-i}-y\opt) + (G-\IM{})\bar{e}_k + e_{k+1}
 	\nonumber\\
	&=
		G\sum_{i=0}^{m_k} \alpha_i^k (G-\IM{})(y_{k-i}-y\opt) + (G-\IM{})\bar{e}_k + e_{k+1}
 	\nonumber\\ 	
 	&=  		
 		G\sum_{i=0}^{m_k}  \alpha_i^k\left[F(y_{k-i}) - e_{k-i}\right]
 		+ (G-\IM{})\bar{e}_k + e_{k+1}
 	\nonumber\\
 	&=
 		G\sum_{i=0}^{m_k}  \alpha_i^k F(y_{k-i}) -\bar{e}_k + e_{k+1}.
\end{align}
By the definition of $\alpha^k$, $\ltwo{\sum_{i=0}^{m_k}  \alpha_i^k F(y_{k-i}) } \leq \ltwo{F(y_k)}$, so \eqref{eq:Fnext:expansion}, \eqref{eq:ynext:term:1}, and \eqref{eq:enext} imply that
\begin{align*}
	\ltwo{F(y_{k+1})} (1- \frac{\varepsilon}{1-\rho}) 
	&\leq
		\rho \ltwo{F(y_k)} 
		+
		\frac{\varepsilon M_\alpha}{1-\rho}\hat{\rho}^{k-m}\ltwo{F(y_0)}
	\\
	&\leq
	\left(
		\frac{\rho}{\hat{\rho}}
		+
		\frac{\varepsilon M_\alpha \,\hat{\rho}^{-m-1}}{1-\rho}
	\right)
	\hat{\rho}^{k+1}\ltwo{F(y_0)}
	\\
	&\leq 	\hat{\rho}^{k+1}\ltwo{F(y_0)},
\end{align*}
where we used the induction hypothesis and \eqref{eq:y:init:cond:r}. 
Appealing to \eqref{eq:Fy:twoside}, the non-expansiveness of $\proxmap_{\gamma h}$, and the fact that $\varepsilon < (1-\ltwo{G})/2$, we obtain
\begin{align*}
	\norm{x_k-x^\star}_2 \leq
	\norm{y_k-y^\star}_2 \leq \frac{1+\rho}{1-\rho}\hat{\rho}^k\norm{y_0-y^{\star}}_2
	\leq 
		\frac{3+\ltwo{G}}{1-\ltwo{G}}\hat{\rho}^k\norm{y_0-y^{\star}}_2,
\end{align*}
which yields the first claim in the theorem.
Finally, the second claim follows by noting that $\varepsilon$ and $\varepsilon'$ are arbitrary and the fact that 
$\lim_{k\to\infty}{a^{1/k}}=1$ for any positive number $a$. This completes the proof.

\section{Proof of Proposition~\ref{prop:example}}\label{appendix:prop1}
We start by recalling the following useful result. For $a, b \in \R$ satisfying $a\neq b$, the solution to the minimization problem
\begin{align*}
	\begin{array}{ll}
		\underset{\alpha_0, \alpha_1 \in\R}{\minimize} & (\alpha_0 a + \alpha_1 b)^2 \\
			\subjectto &  \alpha_0+\alpha_1=1,
	\end{array}
\end{align*}
is given by 
\begin{align}\label{eq:subprob:example:sol}
	\alpha_0 = \frac{b}{b-a} \quad \mbox{and}\quad \alpha_1 = \frac{-a}{b-a}.
\end{align}
Recall also that the AA-GD method is the application of Algorithm~\ref{alg:aa} to the mapping $g(x)=x-\gamma \grad{f(x)}$. Since $m=1$, the $k$-th subproblem ($k\geq 1$) in Step~5 of Algorithm~\ref{alg:aa} boils down to \added{computing} 
\begin{align*}
	\alpha^k = \argmin_{\alpha_0+\alpha_1=1} \left(\alpha_0 \grad{f(x_k)} + \alpha_1 \grad{f(x_{k-1})}\right)^2,
\end{align*}
which together with \eqref{eq:subprob:example:sol} imply that
\begin{align*}
	\alpha_0^k = \frac{\grad{f(x_{k-1})}}{\grad{f(x_{k-1})}- \grad{f(x_k)}}	
	\quad \mbox{and}\quad
	\alpha_1^k = \frac{-\grad{f(x_k)}}{\grad{f(x_{k-1})}- \grad{f(x_k)}}.
\end{align*}
Consequently, we can explicitly compute the next iterate defined in Step~6 of Algorithm~\ref{alg:aa} as
\begin{align}\label{eq:aa:example:iter}
	x_{k+1} 
	= 
		\alpha_0^k g(x_k) + \alpha_1^k g(x_{k-1})
	= 
		\frac{\grad{f(x_{k-1})}}{\grad{f(x_{k-1})}- \grad{f(x_k)}}x_k 
		-
		\frac{\grad{f(x_k)}}{\grad{f(x_{k-1})}- \grad{f(x_k)}}x_{k-1}.
\end{align}

By the construction of $\grad{f(x)}$ and \eqref{eq:aa:example:iter}, it follows that whenever $x_k$ and $x_{k-1}$ belong to the interval $[1, +\infty)$, the next iterate $x_{k+1}$ will take the value $-249$. Similarly, if $x_k$ and $x_{k-1}$ belong to the interval $(-\infty, -1]$, then $x_{k+1}=+249$.  
This motivates us to select the initial interval so that some subsequnece of $\{x_k\}$ will always take the value $+249$ or $-249$, and hence never converge to the origin. To do so, let us examine the pattern of the first few iterates.

First, \replaced{let}{we would like} $x_0, x_1 \in [1, +\infty)$ so that $x_2=-249$.  Given $x_1$ and $x_2$, it is easy to verify that
\begin{align*}
	x_3 = \frac{249(x_1+x_2)}{x_1-x_2+498} = \frac{249(x_1-249)}{x_1+747}.
\end{align*}
Since $x_2 < -1$, if we \replaced{ensure that}{enforce} $x_3\leq-1$, \replaced{we will}{then we must} have $x_4=+249$. Note that for $ x_1 \geq 1$, the right-hand-side of the preceding equation is an increasing function of $x_1$, therefore $x_3< -1$ when $x_1\leq 245$. Also, since $x_1\geq 1$, we have $x_3> -83$. In summary, for $x_1\in [1, 245]$, we have $x_3\in(-83, -1)$ and $x_4=+249$. 
A similar calculation yields
\begin{align*}
x_5 = \frac{-249(x_3+x_4)}{x_3-x_4-498} = \frac{-249(x_3+249)}{x_3-747}.
\end{align*}
Similarly, for $x_3\in(-83, -1)$, $x_5$ is an increasing function of $x_3$, therefore $x_5\in (49.8, 83)$.
Now, since $x_4, x_5 \in [1,\infty)$, $x_6=-249$. The process is now repeated with $x_1$ replaced by $x_5$, $x_2$ replaced by $x_6$, and so on. Note that since $x_5\in (49.8, 83)\subset [1,245]$, all the above results are still valid and can be summarized as:
\begin{align*}
	x_{4n+3} \in (-83,-1), \quad 
	x_{4n+4} =+249, \quad 
	x_{4n+5} \in [1, 245],\quad
	x_{4n+6} = -249, \quad \mbox{for} \quad n=0,1, 2, \ldots,
\end{align*}
which implies that AA-GD will never converges to the optimal solution. 

Indeed, it can be shown that all the four subsequences above will eventually converge.  \replaced{Under our}{Note that under our} initial condition, for  $n=0,1, 2, \ldots$, the iterates $x_{4n+3}$ and $x_{4n+5}$ have the forms
\begin{align*}
	x_{4n+3} &= \frac{249(x_{4n+1}-249)}{x_{4n+1}+747}\\
	x_{4n+5} &= \frac{-249(x_{4n+3}+249)}{x_{4n+3}-747}.
\end{align*}
Thus, we can find a transformation from $x_{4n+1}$ to $x_{4n+5}$ as
\begin{align*}
x_{4n+5} &= \frac{249(x_{4n+1}+249)}{x_{4n+1}+1245}.
\end{align*}
Define $y_n=x_{4n+1}$, then the previous equation can be seen as \replaced{a}{the} fixed-point iteration $y_{n+1}=G(y_n)$ with $G(y):=249(y+249)/(y+1245)$. It is easy to verify that for $y\in[1,245]$, the mapping $G$ is contractive, and hence $\{y_n\}$ converges to the unique fixed-point of $G$ in $[1,245]$, which is $+249(\sqrt{5}-2)$. A parallel argument yields $x_{4n+3}\to -249(\sqrt{5}-2)$  as $n \to \infty$. 

Finally, since $x_1= x_0-(1/L)\grad{f(x_0)}$, to guarantee $x_1\in[1,245]$, a sufficient condition is $x_0\in[2.01, 246.98]$. This completes the proof.

\end{document}